\documentclass[11pt]{article}
\usepackage{mathrsfs}
\usepackage{latexsym,amsmath,amssymb,amsfonts,epsfig,cite,psfrag}
\usepackage{eepic,color,colordvi,amscd,indentfirst,graphics}

\newtheorem{theo}{Theorem}[section]

\newtheorem{lem}[theo]{Lemma}

\newtheorem{prop}[theo]{Proposition}

\newtheorem{claim}[theo]{Claim}
\def\qed{\hfill \rule{4pt}{7pt}}
\def\pf{\noindent {\it Proof.} }

\begin{document}

\title{The saturation number of $W_4$}
\author{Ning Song$^1$\footnote{E-mail:
songning@sdut.edu.cn.},\, Jinze Hu$^1$\footnote{E-mail:
hujinzeh@163.com.},\, Shengjin Ji$^1$\footnote{E-mail:
jishengjin@sdut.edu.cn.}\, and\, Qing Cui$^2$\footnote{E-mail:
cui@nuaa.edu.cn (Corresponding author).}\\
\\
$^1$School of Mathematics and Statistics\\
Shandong University of Technology\\
Zibo 255049, P.R. China\\
\\
$^2$School of Mathematics\\
Nanjing University of Aeronautics and Astronautics\\
Nanjing 210016, P.R. China}
\date{}
\maketitle{}

\begin{abstract}
For a fixed graph $H$, a graph $G$ is called $H$-saturated if $G$
does not contain $H$ as a (not necessarily induced) subgraph, but
$G+e$ contains a copy of $H$ for any $e\in E(\overline{G})$. The
saturation number of $H$, denoted by ${\rm sat}(n,H)$, is the
minimum number of edges in an $n$-vertex $H$-saturated graph. A
wheel $W_n$ is a graph obtained from a cycle of length $n$ by adding
a new vertex and joining it to every vertex of the cycle. A
well-known result of Erd\H{o}s, Hajnal and Moon shows that ${\rm
sat}(n,W_3)=2n-3$ for all $n\geq 4$ and $K_2\vee \overline{K_{n-2}}$
is the unique extremal graph, where $\vee$ denotes the graph join
operation. In this paper, we study the saturation number of $W_4$.
We prove that ${\rm sat}(n,W_4)=\lfloor\frac{5n-10}{2}\rfloor$ for
all $n\geq 6$ and give a complete characterization of the extremal
graphs.
\end{abstract}

\hspace{10pt}{\bf Keywords:} saturation number, wheel, extremal
graph, minimum degree

\hspace{10pt}{\bf Mathematics Subject Classification:} 05C35

\newpage

\section{Introduction}

In this paper we only consider finite simple graphs. For a graph
$G$, we use $V(G)$, $E(G)$, $v(G)$ and $e(G)$ to denote the vertex
set, the edge set, the number of vertices and the number of edges of
$G$, respectively. Let $\overline{G}$ denote the complement graph of
$G$. For any $v\in V(G)$, let $N_G(v)$ and $d_G(v)$ denote the
neighborhood and the degree of $v$ in $G$, respectively, and let
$N_G[v]=N_G(v)\cup \{v\}$. We may omit the subscript $G$ if it is
clear from the context. A vertex $v\in V(G)$ is called a
\emph{universal vertex} of $G$ if $N_G[v]=V(G)$, and the minimum
degree of $G$ is denoted by $\delta(G)$. For any $S\subseteq V(G)$,
we use $G[S]$ to denote the subgraph of $G$ induced by $S$ and
simply write $e(S)$ instead of $e(G[S])$. For any $A,B\subseteq
V(G)$ with $A\cap B=\emptyset$, let $e(A,B)$ denote the number of
edges of $G$ with one endvertex in $A$ and the other endvertex in
$B$. We use $P_n$, $C_n$, $K_n$ and $S_n$ to denote a path, a cycle,
a complete graph and a star with $n$ vertices, respectively. The
\emph{join} of two graphs $G$ and $H$, denoted by $G\vee H$, is the
graph obtained from the disjoint union of $G$ and $H$ by joining
each vertex of $G$ to each vertex of $H$. A \emph{wheel} $W_n$ is a
graph obtained from a cycle $C_n$ by adding a new vertex $v$ and
joining it to every vertex of $C_n$ (i.e. $W_n=C_n\vee\{v\}$), where
the cycle $C_n$ and the vertex $v$ are called the \emph{rim} and the
\emph{center} of $W_n$, respectively. For any positive integer $k$,
let $[k]$ denote the set $\{1,2,\ldots,k\}$. We write $A:=B$ to
rename $B$ as $A$.

For a fixed graph $H$, a graph is \emph{$H$-free} if it does not
contain $H$ as a (not necessarily induced) subgraph. A graph $G$ is
called \emph{$H$-saturated} if $G$ is $H$-free but $G+e$ contains a
copy of $H$ for any $e\in E(\overline{G})$. The \emph{saturation
number} of $H$, denoted by ${\rm sat}(n,H)$, is the minimum number
of edges in an $n$-vertex $H$-saturated graph. An $n$-vertex
$H$-saturated graph with ${\rm sat}(n,H)$ edges is called an
\emph{extremal} graph for $H$, and the set of all $n$-vertex
extremal graphs for $H$ is denoted by ${\rm Sat}(n,H)$.

The study of the saturation numbers of graphs was initiated by
Erd\H{o}s, Hajnal and Moon in~\cite{EHM64}, in which the authors
proved that ${\rm sat}(n,K_k)=(k-2)n-{k-1\choose 2}$ and
$K_{k-2}\vee \overline{K_{n-k+2}}$ is the unique extremal graph.
Later, K\'{a}szonyi and Tuza~\cite{KT86} showed that ${\rm
sat}(n,H)=O(n)$ for any graph $H$ and determined the exact values of
${\rm sat}(n,H)$ for $H\in\{S_k,P_k,tK_2\}$. Since then, there has
been a large quantity of work in determining the saturation numbers
of various classes of graphs such as
cliques~\cite{CP22,CY24,FFGJ09b},
cycles~\cite{Ch09,Ch11,FJMTW12,FK13,GKK24,LSWZ21,MHHG21,Ol72,Tu89},
complete multipartite
graphs~\cite{CFG08,Ch13,GS07,HL21,HLSZ24,SW17},
trees~\cite{DW04,FFGJ09a} and
forests~\cite{CLLYZ23,CFFGJM15,FW15,HLL23,LHL23}. However, the exact
value of ${\rm sat}(n,H)$ and a complete characterization of ${\rm
Sat}(n,H)$ are known for very few special classes of graphs $H$. We
refer the readers to the nice survey of Currie, Faudree, Faudree and
Schmitt~\cite{CFFS21} for a summary of known results on saturation
numbers.

In this paper, we are interested in studying the saturation numbers
of wheels. Notice that $W_3=K_4$, the aforementioned result of
Erd\H{o}s, Hajnal and Moon~\cite{EHM64} implies that ${\rm
sat}(n,W_3)=2n-3$ for all $n\geq 4$ and $K_2\vee \overline{K_{n-2}}$
is the unique extremal graph. As far as we are aware, this is the
only known result for wheels so far. As a natural next step, the aim
of this paper is to determine the exact value of ${\rm sat}(n,W_4)$
and give a complete characterization of ${\rm Sat}(n,W_4)$  for all
$n\geq 5$.

\begin{figure}[ht]
\begin{center}
\includegraphics[scale=0.8]{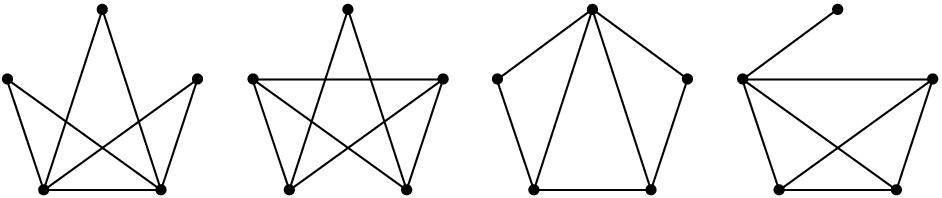}
\caption{The four graphs with $5$ vertices and $7$ edges.}
\label{fig1}
\end{center}
\end{figure}

\begin{figure}[ht]
\begin{center}
\includegraphics[scale=0.8]{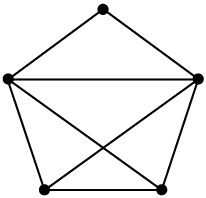}
\caption{The extremal graph $H^*$.}
\label{fig2}
\end{center}
\end{figure}

We point out here that an easy argument can show that ${\rm
sat}(n,W_4)=8$ when $n=5$ and the extremal graph is unique. Since
$e(W_4)=8$, we know that every $W_4$-saturated graph contains at
least $7$ edges. Note that there are exactly four graphs with $5$
vertices and $7$ edges, none of which is $W_4$-saturated (see
Figure~\ref{fig1}). On the other hand, $W_4$ and $H^*$ are the only
two graphs with $5$ vertices and $8$ edges, where $H^*$ is the graph
obtained from $K_5$ by deleting two consecutive edges (see
Figure~\ref{fig2}). Since $H^*$ is $W_4$-saturated, we conclude that
${\rm sat}(5,W_4)=8$ and ${\rm Sat}(5,W_4)=\{H^*\}$.

Hence, we need only to consider $n\geq 6$ in the following
arguments. In order to state our main result, we need to introduce
several families of graphs.

\vspace{15pt}
\begin{figure}[ht]
\begin{center}
\includegraphics[scale=0.8]{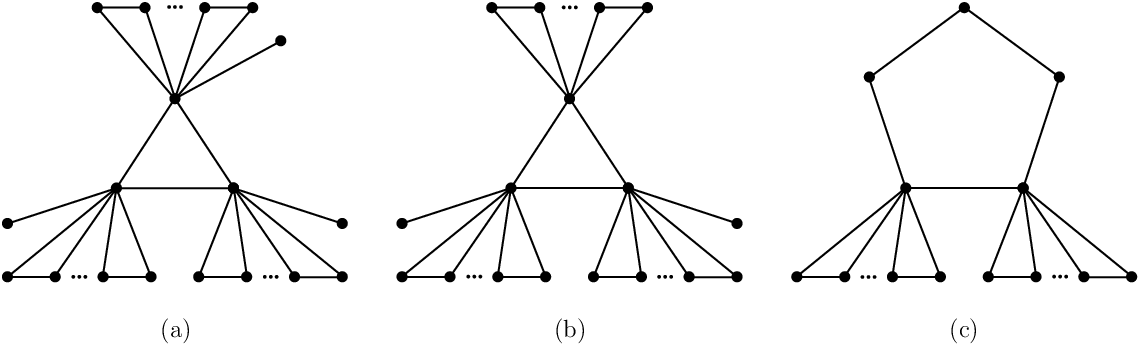}
\caption{The graph families $\mathcal{F}_n^1$, $\mathcal{F}_n^2$ and
$\mathcal{F}_n^3$.} \label{fig3}
\end{center}
\end{figure}

For any even integer $n\geq 6$, let $\mathcal{F}_n^1$ denote the
family of $n$-vertex graphs $F$ such that $F$ has a `central'
triangle, each of whose vertices is adjacent to exactly one vertex
of degree $1$, and the remaining $n-6$ vertices of $F$ are in
adjacent pairs, each of them joined to a vertex of the central
triangle (see Figure~\ref{fig3}(a) for an illustration). For any odd
integer $n\geq 5$, let $\mathcal{F}_n^2$ denote the family of
$n$-vertex graphs which are obtained from the graphs in
$\mathcal{F}_{n+1}^1$ by deleting one vertex of degree $1$ (see
Figure~\ref{fig3}(b)). For any odd integer $n\geq 5$, let
$\mathcal{F}_n^3$ denote the family of $n$-vertex graphs $F$ such
that $F$ consists of a $C_5$, two consecutive vertices of which are
joined to arbitrary numbers of adjacent pairs (see
Figure~\ref{fig3}(c) for an illustration). These families of graphs
were first introduced by Ollmann in~\cite{Ol72}, in which the author
determined ${\rm sat}(n,C_4)$ and ${\rm Sat}(n,C_4)$ for all $n\geq
5$. (An alternative proof was later given by Tuza~\cite{Tu89}.)

\begin{theo}{\rm (Ollmann~\cite{Ol72}, Tuza~\cite{Tu89})}
\label{theo1.1}
For $n\geq 5$,
${\rm sat}(n,C_4)=\lfloor\frac{3n-5}{2}\rfloor$ and
\begin{align*}
{\rm Sat}(n,C_4)=\left\{\begin{array}{ll}\mathcal{F}_n^1,
& \;\mbox{ if } n \mbox{ is even},\\
\mathcal{F}_n^2\cup\mathcal{F}_n^3, & \;\mbox{ if } n  \mbox{ is odd}.
\end{array}\right.
\end{align*}
\end{theo}

For any odd integer $n\geq 7$, we define $\mathcal{A}_n^1:=\{F\vee
K_1: F\in\mathcal{F}_{n-1}^1\}$. For any even integer $n\geq 6$, we
define $\mathcal{A}_n^2:=\{F\vee K_1: F\in\mathcal{F}_{n-1}^2\}$ and
$\mathcal{A}_n^3:=\{F\vee K_1: F\in\mathcal{F}_{n-1}^3\}$.

\vspace{15pt}
\begin{figure}[ht]
\begin{center}
\includegraphics[scale=0.85]{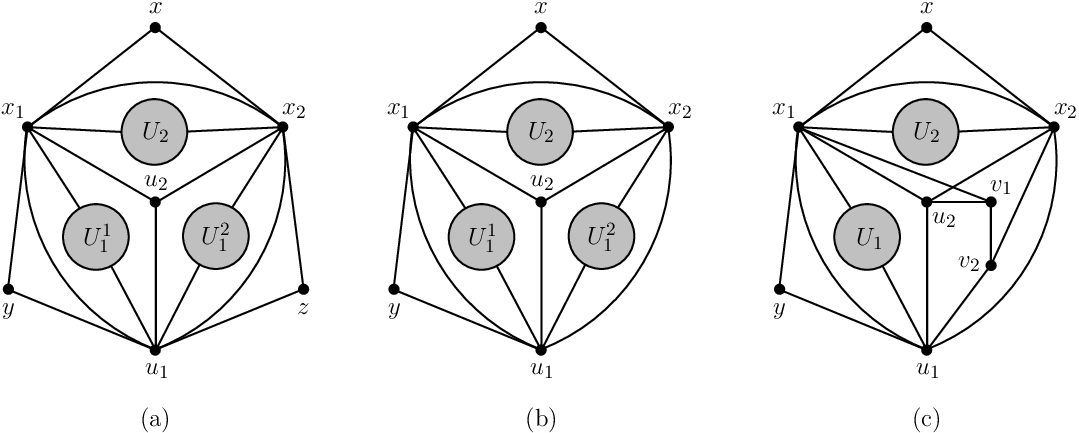}
\caption{The graph families $\mathcal{B}_n^1$, $\mathcal{B}_n^2$ and
$\mathcal{B}_n^3$.}
\label{fig4}
\end{center}
\end{figure}

For any odd integer $n\geq 7$, let $\mathcal{B}_n^1$ denote the
family of $n$-vertex graphs $G$ such that
$V(G)=\{x,y,z,x_1,x_2,u_1,u_2\}\cup U_1^1\cup U_1^2\cup U_2$ (it is
possible that $U_1^1$, $U_1^2$ or $U_2$ is empty) and the
following properties hold:
\begin{itemize}
\item [(i)] $G[\{x_1,x_2,u_1,u_2\}]\cong K_4$;
\item [(ii)] $d(x)=d(y)=d(z)=2$ with $N(x)=\{x_1,x_2\}$,
$N(y)=\{x_1,u_1\}$ and $N(z)=\{x_2,u_1\}$;
\item [(iii)] $G[U]$ is a matching for any $U\in\{U_1^1, U_1^2, U_2\}$ with $U\neq\emptyset$;
\item [(iv)] every vertex in $U_1^i$ is adjacent to both $x_i$ and $u_1$ for each $i\in [2]$
with $U_1^i\neq\emptyset$;
\item [(v)] every vertex in $U_2$ is adjacent to both $x_1$ and $x_2$ if $U_2\neq\emptyset$.
\end{itemize}
See Figure~\ref{fig4}(a) for an illustration. For any even integer
$n\geq 8$, let $\mathcal{B}_n^2$ denote the family of $n$-vertex
graphs which are obtained from the graphs in $\mathcal{B}_{n+1}^1$
with $U_1^2\neq\emptyset$ by deleting the vertex $z$ (see
Figure~\ref{fig4}(b)). For any even integer $n\geq 8$, let
$\mathcal{B}_n^3$ denote the family of $n$-vertex graphs $G$ such
that $V(G)=\{x,y,x_1,x_2,y_1,y_2,u_1,u_2\}\cup U_1\cup U_2$ (it is
possible that $U_1$ or $U_2$ is empty) and the following properties
hold:
\begin{itemize}
\item [(i)] $G[\{x_1,x_2,u_1,u_2\}]\cong K_4$;
\item [(ii)] $d(x)=d(y)=2$ with $N(x)=\{x_1,x_2\}$ and $N(y)=\{x_1,u_1\}$;
\item [(iii)] $d(v_1)=d(v_2)=3$ with $N(v_1)=\{x_1,v_2,u_2\}$ and $N(v_2)=\{x_2,v_1,u_1\}$;
\item [(iv)] $G[U]$ is a matching for any $U\in\{U_1, U_2\}$ with $U\neq\emptyset$;
\item [(v)] every vertex in $U_1$ is adjacent to both $x_1$ and $u_1$ if $U_1\neq\emptyset$;
\item [(vi)] every vertex in $U_2$ is adjacent to both $x_1$ and $x_2$ if $U_2\neq\emptyset$.
\end{itemize}
Please refer to Figure~\ref{fig4}(c) for a detailed illustration.

We can now state the main result of this paper.

\begin{theo}\label{theo1.2}
For $n\geq 6$,
${\rm sat}(n,W_4)=\lfloor\frac{5n-10}{2}\rfloor$ and
\begin{align*}
{\rm
Sat}(n,W_4)=\left\{\begin{array}{ll}\mathcal{A}_n^1\cup\mathcal{B}_n^1,
& \;\mbox{ if } n \mbox{ is odd},\\
\mathcal{A}_n^2\cup\mathcal{A}_n^3\cup\mathcal{B}_n^2\cup\mathcal{B}_n^3,
& \;\mbox{ if } n  \mbox{ is even}.
\end{array}\right.
\end{align*}
\end{theo}

The rest of this paper is organized as follows. In Section 2, we
show that all graphs in $\mathcal{A}_n^1$, $\mathcal{A}_n^2$,
$\mathcal{A}_n^3$, $\mathcal{B}_n^1$, $\mathcal{B}_n^2$ and
$\mathcal{B}_n^3$ are $W_4$-saturated. In Section 3, we investigate
some properties of $W_4$-saturated graphs. The proof of
Theorem~\ref{theo1.2} will be given in Section 4.

\section{The upper bound}

In this section, we shall prove that all graphs in
$\mathcal{A}_n^1$, $\mathcal{A}_n^2$, $\mathcal{A}_n^3$,
$\mathcal{B}_n^1$, $\mathcal{B}_n^2$ and $\mathcal{B}_n^3$ are
$W_4$-saturated and contain exactly $\lfloor\frac{5n-10}{2}\rfloor$
edges, which implies that $\lfloor\frac{5n-10}{2}\rfloor$ is an
upper bound of ${\rm sat}(n,W_4)$.

\begin{prop}\label{prop2.1}
For any odd integer $n\geq 7$, the graphs in $\mathcal{A}_n^1$ are
$W_4$-saturated and contain $\lfloor\frac{5n-10}{2}\rfloor$ edges.
\end{prop}

\pf Let $G$ be a graph in $\mathcal{A}_n^1$. Then by the definition
of $\mathcal{A}_n^1$, we may assume that $G=F\vee \{v\}$ for some
$F\in\mathcal{F}_{n-1}^1$ and $v$ is a universal vertex of $G$. By
Theorem~\ref{theo1.1}, we can know that $F$ is $C_4$-saturated and
$e(F)=\lfloor\frac{3(n-1)-5}{2}\rfloor=\lfloor\frac{3n-8}{2}\rfloor$.
Hence, we have
\begin{align*}
e(G)=e(F)+(n-1)=\lfloor\frac{3n-8}{2}\rfloor+(n-1)=\lfloor\frac{5n-10}{2}\rfloor.
\end{align*}

Next, we show that $G$ is $W_4$-free. Suppose not, and let $H$ be a
copy of $W_4$ of $G$. Since $F$ is $C_4$-free, we see that $F$ is
also $W_4$-free and thus $v\in V(H)$. Observe that $H-v$ contains a
copy of $C_4$ (no matter whether $v$ is the center of $H$ or not)
and $H-v\subseteq F$, we derive a contradiction to the fact that $F$
is $C_4$-free. Therefore, $G$ is $W_4$-free.

Finally, we show that $G$ is $W_4$-saturated. Let $st$ be an edge in
$\overline{G}$. Then we have $s,t\in V(F)$ (because $v$ is a
universal vertex of $G$). Since $F$ is $C_4$-saturated, there exists
a copy of $C_4$ in $F+st$, say $R$. Then the subgraph of $G+st$
induced by $V(R)\cup\{v\}$ contains a copy of $W_4$. Thus, $G$ is
$W_4$-saturated. \qed

\begin{prop}\label{prop2.2}
For any even integer $n\geq 6$, the graphs in
$\mathcal{A}_n^2\cup\mathcal{A}_n^3$ are $W_4$-saturated and contain
$\lfloor\frac{5n-10}{2}\rfloor$ edges.
\end{prop}

\pf The proof is the same as that of Proposition~\ref{prop2.1}.
\qed

\begin{prop}\label{prop2.3}
For any odd integer $n\geq 7$, the graphs in $\mathcal{B}_n^1$ are
$W_4$-saturated and contain $\lfloor\frac{5n-10}{2}\rfloor$ edges.
\end{prop}

\pf Let $G$ be a graph in $\mathcal{B}_n^1$, where the vertices of
$G$ are labeled as shown in Figure~\ref{fig4}(a). Since $n$ is odd,
it follows from the definition of $\mathcal{B}_n^1$ that
\begin{align*}
e(G)=6+6+\frac{n-7}{2}+(n-7)\cdot
2=\frac{5n-11}{2}=\lfloor\frac{5n-10}{2}\rfloor.
\end{align*}

Next, we show that $G$ is $W_4$-free. Suppose not. Let $H$ be a copy
of $W_4$ of $G$ and let $R$ be the rim of $H$. Notice that $x_1$,
$x_2$ and $u_1$ are the only three vertices of $G$ with degree at
least $4$, we may assume by symmetry that $u_1$ is the center of
$H$. Then $V(R)\subseteq N(u_1)$. Moreover, since both $x_1$ and
$x_2$ are cut-vertices of $G[N(u_1)]$, we can further conclude that
either $V(R)\subseteq U_1^1\cup\{x_1\}$ or $V(R)\subseteq
U_1^2\cup\{x_2\}$. But this is impossible since it is easy to
observe that neither $G[U_1^1\cup\{x_1\}]$ nor $G[U_1^2\cup\{x_2\}]$
contains a copy of $C_4$, a contradiction. Hence, $G$ is $W_4$-free.

Finally, we show that $G$ is $W_4$-saturated. Let $st$ be an edge in
$\overline{G}$ and let $G':=G+st$. By symmetry, we need only to
consider the following cases.
\begin{itemize}
\item [(i)] If $s=x$ and $t\in\{y,z,u_2\}\cup U_1^1\cup U_1^2$,
then $G'[\{x_1,x_2,u_1,s,t\}]$ contains a copy of $W_4$.
\item [(ii)] If $s=x$ and $t=u_1$, then $G'[\{x_1,x_2,u_2,s,t\}]$
contains a copy of $W_4$.
\item [(iii)] If $s\in\{x,u_1\}$ and $t\in U_2$, then $G'[\{x_1,x_2,s,t,t'\}]$
contains a copy of $W_4$, where $t'$ is the unique neighbor of $t$
in $U_2$.
\item [(iv)] If $s=u_2$ and $t\in U_1^1\cup U_1^2\cup U_2$,
then $G'[\{x_1,x_2,u_1,s,t\}]$ contains a copy of $W_4$.
\item [(v)] If $s\in U_2$ and
$t\in U_1^1\cup U_1^2$, then $G'[\{x_1,x_2,u_1,s,t\}]$ contains a
copy of $W_4$.
\item [(vi)] If $s,t\in U_2$, then $G'[\{x_1,x_2,s,t,t'\}]$
contains a copy of $W_4$, where $t'$ is the unique neighbor of $t$
in $U_2$.
\end{itemize}
In all cases, we see that $G'$ contains a copy of $W_4$. Thus, $G$
is $W_4$-saturated. \qed

\begin{prop}\label{prop2.4}
For any even integer $n\geq 8$, the graphs in $\mathcal{B}_n^2$ are
$W_4$-saturated and contain $\lfloor\frac{5n-10}{2}\rfloor$ edges.
\end{prop}

\pf Let $G$ be a graph in $\mathcal{B}_n^2$, where the vertices of
$G$ are labeled as shown in Figure~\ref{fig4}(b). Then by the
definition of $\mathcal{B}_n^2$, we may assume that $G=F-z$ for some
$F\in\mathcal{B}_{n+1}^1$ with $U_1^2\neq\emptyset$. Since $n$ is
even and by Proposition~\ref{prop2.3}, we know that $F$ is
$W_4$-saturated and
$e(F)=\lfloor\frac{5(n+1)-10}{2}\rfloor=\lfloor\frac{5n-5}{2}\rfloor=\lfloor\frac{5n-6}{2}\rfloor$.
This implies that $G$ is $W_4$-free (because $F$ is $W_4$-free and
$G\subseteq F$) and
\begin{align*}
e(G)=e(F)-2=\lfloor\frac{5n-6}{2}\rfloor-2=\lfloor\frac{5n-10}{2}\rfloor.
\end{align*}

Let $st$ be an edge in $\overline{G}$. Then $st\in E(\overline{F})$
and it follows from $F$ is $W_4$-saturated that there exists a copy
of $W_4$ in $F+st$, say $H$. Since $d_F(z)=2$, we have $z\notin
V(H)$. This means that $H$ is also a subgraph of $G+st$. Therefore,
$G$ is $W_4$-saturated.   \qed

\begin{prop}\label{prop2.5}
For any even integer $n\geq 8$, the graphs in $\mathcal{B}_n^3$ are
$W_4$-saturated and contain $\lfloor\frac{5n-10}{2}\rfloor$ edges.
\end{prop}

\pf Let $G$ be a graph in $\mathcal{B}_n^3$, where the vertices of
$G$ are labeled as shown in Figure~\ref{fig4}(c). By the definition
of $\mathcal{B}_n^3$, we derive that
\begin{align*}
e(G)=6+4+5+\frac{n-8}{2}+(n-8)\cdot
2=\frac{5n-10}{2}=\lfloor\frac{5n-10}{2}\rfloor.
\end{align*}

Next, we show that $G$ is $W_4$-free. Suppose not. Let $H$ be a copy
of $W_4$ of $G$ and let $R$ be the rim of $H$. Note that $x_1$,
$x_2$, $u_1$ and $u_2$ are the only four vertices of $G$ with degree
at least $4$. Since $d(u_2)=4$ and $G[N(u_2)]$ does not contain a
copy of $C_4$, we conclude that $u_2$ is not the center of $H$.
Suppose $x_1$ is the center of $H$. Then $V(R)\subseteq N(x_1)$.
Since both $u_1$ and $x_2$ are cut-vertices of $G[N(x_1)]$, we
deduce that either $V(R)\subseteq U_1\cup\{u_1\}$ or $V(R)\subseteq
U_2\cup\{x_2\}$. However, this is impossible since it is easy to see
that neither $G[U_1\cup\{u_1\}]$ nor $G[U_2\cup\{x_2\}]$ contains a
copy of $C_4$, a contradiction. Similarly, we can show that neither
$x_2$ nor $u_1$ is the center of $H$. Thus, $G$ is $W_4$-free.

Finally, we show that $G$ is $W_4$-saturated. Let $st$ be an edge in
$\overline{G}$ and let $G':=G+st$. By symmetry, it suffices to
consider the following cases.
\begin{itemize}
\item [(i)] If $s=x$ and $t\in\{y,u_2,v_2\}\cup U_1$,
then $G'[\{x_1,x_2,u_1,s,t\}]$ contains a copy of $W_4$.
\item [(ii)] If $s=x$ and $t\in\{u_1,v_1\}$, then $G'[\{x_1,x_2,u_2,s,t\}]$
contains a copy of $W_4$.
\item [(iii)] If $s\in\{x,u_1\}$ and $t\in U_2$, then $G'[\{x_1,x_2,s,t,t'\}]$
contains a copy of $W_4$, where $t'$ is the unique neighbor of $t$
in $U_2$.
\item [(iv)] If $s=u_1$ and
$t=v_1$, then $G'[\{x_1,x_2,u_2,s,t\}]$ contains a copy of $W_4$.
\item [(v)] If $s=u_2$ and $t\in \{v_2\}\cup U_1\cup U_2$,
then $G'[\{x_1,x_2,u_1,s,t\}]$ contains a copy of $W_4$.
\item [(vi)] If $s=v_1$ and
$t\in U_1\cup U_2$, then $G'[\{x_1,x_2,u_2,s,t\}]$ contains a copy
of $W_4$.
\item [(vii)] If $s=v_2$ and $t\in U_1\cup U_2$,
then $G'[\{x_1,x_2,u_1,s,t\}]$ contains a copy of $W_4$.
\item [(viii)] If $s=v_2$ and $t=x_1$,
then $G'[\{x_2,u_1,u_2,s,t\}]$ contains a copy of $W_4$.
\item [(ix)] If $s\in U_2$ and
$t\in U_1$, then $G'[\{x_1,x_2,u_1,s,t\}]$ contains a copy of $W_4$.
\item [(x)] If $s,t\in U_2$, then $G'[\{x_1,x_2,s,t,t'\}]$
contains a copy of $W_4$, where $t'$ is the unique neighbor of $t$
in $U_2$.
\end{itemize}
In all cases, we see that $G'$ contains a copy of $W_4$. Therefore,
$G$ is $W_4$-saturated. \qed

\section{Properties of $W_4$-saturated graphs}

In this section, we investigate some useful properties of
$W_4$-saturated graphs and define two functions on the set of
vertices of $W_4$-saturated graphs. These will be used in the next
section to prove the main result of this paper.

Fix a $W_4$-saturated graph $G$ with $n\geq 6$ vertices. Clearly,
$G\not\cong K_n$. We choose a vertex $x$ in $G$ such that
$d(x)=\delta(G)$ and $e(N[x])$ is as small as possible. Let
$N(x)=\{x_1,x_2,\dots,x_{\delta(G)}\}$ and $V_x:=V(G)\setminus
N[x]$. Then $V_x\neq\emptyset$. For each $i=0,1,\ldots,\delta(G)$,
we define $V_i:=\{v\in V_x:|N(v)\cap N(x)|=i\}$.

\begin{lem}\label{lem3.1}
The following statements hold:
\begin{itemize}
\item [$($i$)$] $\delta(G)\geq 2$;
\item [$($ii$)$] for any pair of non-adjacent vertices $s$ and $t$ in $G$, we have
$N(s)\cap N(t)\neq\emptyset$ (i.e. $s$ and $t$ have at least one
common neighbor);
\item [$($iii$)$] $V_0=\emptyset$.
\end{itemize}
\end{lem}

\pf Let $v$ be a vertex in $V_x$. Since $G$ is $W_4$-saturated,
there exists a copy of $W_4$ in $G+vx$, say $H$. It is clear that
$vx\in E(H)$ and $3\leq d_H(x)\leq 4$. Since $H-vx\subseteq G$, we
know that $\delta(G)=d_G(x)\geq d_H(x)-1\geq 2$. So we have (i).

Suppose $s$ and $t$ are two non-adjacent vertices in $G$. Let $H'$
be a copy of $W_4$ in $G+st$. It is easy to observe that $s$ and $t$
have at least one common neighbor in $H'$ (no matter whether $s$ or
$t$ is the center of $H'$ or not). Since $H'-st\subseteq G$, we can
derive that any common neighbor of $s$ and $t$ in $H'$ is also a
common neighbor of them in $G$. Hence, $N(s)\cap N(t)\neq\emptyset$.
This proves (ii).

It follows from (ii) that $N(v)\cap N(x)\neq\emptyset$ for any $v\in
V_x$. Thus, $V_0=\emptyset$. This proves (iii). \qed

\vspace{15pt}
\begin{figure}[ht]
\begin{center}
\includegraphics[scale=0.8]{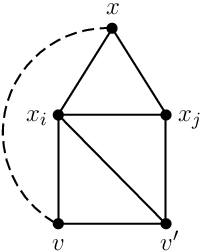}
\caption{The configuration in Lemma~\ref{lem3.2}.}
\label{fig5}
\end{center}
\end{figure}
\vspace{-10pt}

\begin{lem}\label{lem3.2}
Let $v$ be a vertex in $V_1$ such that $N(v)\cap N(x)=\{x_i\}$ for
some $i\in [\delta(G)]$. Then there exist some $x_j$ ($j\neq i$) and
$v'\in V_2\cup\cdots\cup V_{\delta(G)}$ such that
$vv',v'x_i,v'x_j,x_ix_j\in E(G)$ (see Figure~\ref{fig5}).
\end{lem}

\pf Since $G$ is $W_4$-saturated, there exists a copy of $W_4$ in
$G+vx$, say $H$. Notice that $x_i$ is the unique common neighbor of
$v$ and $x$ in $G$, we conclude that the center of $H$ must be
$x_i$. Let $vxx'v'v$ be the rim of $H$. Then
$xx',vv',v'x_i,v'x',x_ix'\in E(G)$. Since $xx',x_ix'\in E(G)$, we
have $x'=x_j$ for some $j\in [\delta(G)]\setminus\{i\}$. Moreover,
because $v\in V_1$ and $vv',vx_i,v'x_i,v'x_j\in E(G)$, we deduce
that $v'\in V_2\cup\cdots\cup V_{\delta(G)}$.    \qed

\bigskip

For the sake of brevity, the vertex $v'$ in Lemma~\ref{lem3.2} is
called a \emph{shadow} of $v$. Then Lemma~\ref{lem3.2} shows that
every vertex in $V_1$ has at least one shadow in $V_2\cup\cdots\cup
V_{\delta(G)}$.

\begin{lem}\label{lem3.3}
Suppose $\delta(G)=3$. Let $u$ be a vertex in $V_2$ such that
$N(u)\cap N(x)=\{x_i,x_j\}$ for some $i,j\in [3]$ and $N(u)\cap
(V_2\cup V_3)=\emptyset$. Then $x_kx_i,x_kx_j\in E(G)$, where $k\in
[3]\setminus\{i,j\}$.
\end{lem}

\pf Without loss of generality, we may assume that $\{i,j\}=\{1,2\}$
and $k=3$. Therefore, to prove this lemma, it suffices to show that
$x_3x_1,x_3x_2\in E(G)$. Let $H$ be a copy of $W_4$ in $G+ux$. Since
$N(u)\cap N(x)=\{x_1,x_2\}$, we know that the center of $H$ must be
one vertex in $\{u,x,x_1,x_2\}$.

First, suppose $x$ is the center of $H$. Since $N(u)\cap
N(x)=\{x_1,x_2\}$, we see that the rim of $H$ must be $ux_1x_3x_2u$
and hence $x_3x_1,x_3x_2\in E(G)$, as desired.

Next, suppose $u$ is the center of $H$. Since $N(u)\cap
N(x)=\{x_1,x_2\}$, we have $x_1,x_2\in V(H)$ and $x_3\notin V(H)$.
Let $xx_1vx_2x$ be the rim of $H$. Then $vu,vx_1,vx_2\in E(G)$, and
thus $v\in V_2\cup V_3$. But this contradicts the assumption that
$N(u)\cap (V_2\cup V_3)=\emptyset$.

Finally, suppose by symmetry that $x_1$ is the center of $H$. Let
$uxvwu$ be the rim of $H$. Then $uw,wv,vx,wx_1,vx_1\in E(G)$. Since
$vx,vx_1\in E(G)$, we derive that $v\in\{x_2,x_3\}$. On the other
hand, because $uw\in E(G)$ and $N(u)\cap (V_2\cup V_3)=\emptyset$,
we conclude that $w\notin V_2\cup V_3$. This, together with
$wv,wx_1\in E(G)$, implies that $w\in\{x_2,x_3\}$. Since $N(u)\cap
N(x)=\{x_1,x_2\}$ and $uw\in E(G)$, we have $w=x_2$ and $v=x_3$.
Hence, $x_3x_1,x_3x_2\in E(G)$.     \qed

\bigskip

In the rest of this section, we define two functions which will be
frequently used in Section 4 to give the lower bound of $e(G)$.

The first function is defined as follows: For each $i\in[\delta(G)]$
and each $v\in V_i$, let
\begin{align*}
f(v)=i+0.5|N(v)\cap V_i|+|N(v)\cap (V_{i+1}\cup\cdots\cup
V_{\delta(G)})|. \tag{$1$}
\end{align*}

\begin{lem}\label{lem3.4}
$e(G)=e(N[x])+\sum\limits_{v\in V_x}f(v)$.
\end{lem}

\pf Note that $V_0=\emptyset$ (by Lemma~\ref{lem3.1}(iii)). By the
definition of $f$-function, we know that
\begin{align*}
e(G)&=e(N[x])+e(V_x,N(x))+e(V_x)\\
&=e(N[x])+\sum_{i=1}^{\delta(G)}e(V_i,N(x))+\sum_{i=1}^{\delta(G)}(e(V_i)+e(V_i,V_{i+1}\cup\cdots\cup
V_{\delta(G)}))\\
&=e(N[x])+\sum_{i=1}^{\delta(G)}\sum_{v\in V_i}\left(|N(v)\cap N(x)|
+0.5|N(v)\cap V_i|+|N(v)\cap (V_{i+1}\cup\cdots\cup V_{\delta(G)})|\right)\\
&=e(N[x])+\sum_{i=1}^{\delta(G)}\sum_{v\in V_i}\left(i+0.5|N(v)\cap V_i|
+|N(v)\cap (V_{i+1}\cup\cdots\cup V_{\delta(G)})|\right)\\
&=e(N[x])+\sum_{i=1}^{\delta(G)}\sum_{v\in V_i}f(v)\\
&=e(N[x])+\sum_{v\in V_x}f(v).
\end{align*}
This completes the proof of Lemma~\ref{lem3.4}.    \qed

\bigskip

The second function is defined as follows: For each $v\in V_x$, let
\begin{align*}
g(v)=|N(v)\cap N(x)|+0.5|N(v)\cap V_x|. \tag{$2$}
\end{align*}

\begin{lem}\label{lem3.5}
$e(G)=e(N[x])+\sum\limits_{v\in V_x}g(v)$.
\end{lem}

\pf By the definition of $g$-function, we have
\begin{align*}
e(G)&=e(N[x])+e(V_x,N(x))+e(V_x)\\
&=e(N[x])+\sum_{v\in V_x}|N(v)\cap N(x)|+\sum_{v\in
V_x}0.5|N(v)\cap V_x|\\
&=e(N[x])+\sum_{v\in V_x}\left(|N(v)\cap N(x)|+0.5|N(v)\cap V_x|\right)\\
&=e(N[x])+\sum_{v\in V_x}g(v).
\end{align*}
This proves Lemma~\ref{lem3.5}.    \qed

\section{Proof of Theorem~\ref{theo1.2}}

It follows from
Propositions~\ref{prop2.1},~\ref{prop2.2},~\ref{prop2.3},~\ref{prop2.4}
and~\ref{prop2.5} that $\lfloor\frac{5n-10}{2}\rfloor$ is an upper
bound of ${\rm sat}(n,W_4)$. In the rest of the paper, we shall show
that $\lfloor\frac{5n-10}{2}\rfloor$ is also a lower bound of ${\rm
sat}(n,W_4)$ and characterize the extremal graphs.

Let $G$ be a $W_4$-saturated graph with $n\geq 6$ vertices. In order
to prove Theorem~\ref{theo1.2}, it suffices to show that
$e(G)\geq\lfloor\frac{5n-10}{2}\rfloor$ with equality if and only if
$G\in\mathcal{A}_n^1\cup\mathcal{B}_n^1$ when $n$ is odd and
$G\in\mathcal{A}_n^2\cup\mathcal{A}_n^3\cup\mathcal{B}_n^2\cup\mathcal{B}_n^3$
when $n$ is even. Moreover, since $e(G)$ is an integer, it is easy
to check that $e(G)\geq\lfloor\frac{5n-10}{2}\rfloor$ if and only if
$e(G)\geq\frac{5n-11}{2}$.

Suppose $G$ contains a universal vertex, say $v$. Then $G=F\vee
\{v\}$, where $F$ is an $(n-1)$-vertex graph. Since $G$ is
$W_4$-free, we see that $F$ is $C_4$-free. (Otherwise, suppose $R$
is a copy of $C_4$ of $F$, then $G[V(R)\cup\{v\}]$ contains a copy
of $W_4$, a contradiction.) Let $st$ be an edge in $\overline{F}$
(also in $\overline{G}$). Since $G$ is $W_4$-saturated, there exists
a copy of $W_4$ in $G+st$, say $H$. Note that if $v\notin V(H)$ then
$H\subseteq F+st$ and $H$ contains a copy of $C_4$, and if $v\in
V(H)$ then $H-v\subseteq F+st$ and $H-v$ contains a copy of $C_4$.
In both cases, we can find a copy of $C_4$ in $F+st$. Thus, $F$ is
$C_4$-saturated. Then by Theorem~\ref{theo1.1}, we derive that
$e(F)\geq
\lfloor\frac{3(n-1)-5}{2}\rfloor=\lfloor\frac{3n-8}{2}\rfloor$ with
equality if and only if $F\in\mathcal{F}_{n-1}^1$ when $n-1$ is even
and $F\in\mathcal{F}_{n-1}^2\cup\mathcal{F}_{n-1}^3$ when $n-1$ is
odd. This implies that
\begin{align*}
e(G)=e(F)+(n-1)\geq\lfloor\frac{3n-8}{2}\rfloor+(n-1)=\lfloor\frac{5n-10}{2}\rfloor
\end{align*}
with equality if and only if $G\in\mathcal{A}_n^1$ when $n$ is odd
and $G\in\mathcal{A}_n^2\cup\mathcal{A}_n^3$ when $n$ is even.

Therefore, we may assume that $G$ contains no universal vertex. If
$\delta(G)\geq 5$, then $e(G)\geq
\frac{5n}{2}>\lfloor\frac{5n-10}{2}\rfloor$. Hence by
Lemma~\ref{lem3.1}(i), we may further assume that
$2\leq\delta(G)\leq 4$. Let $x$ be a vertex in $G$ such that
$d(x)=\delta(G)$ and $e(N[x])$ is as small as possible. Let
$N(x)=\{x_1,x_2,\dots,x_{\delta(G)}\}$ and $V_x:=V(G)\setminus
N[x]$. For each $i=0,1,\ldots,\delta(G)$, we define $V_i:=\{v\in
V_x:|N(v)\cap N(x)|=i\}$. Then by Lemma~\ref{lem3.1}(iii), we deduce
that $V_0=\emptyset$.

In the following, we divide the rest of the proof into six parts
according to the values of $\delta(G)$ and $e(N[x])$.

\subsection{$\delta(G)=2$}

\medskip

In this part, $V_2\neq\emptyset$. (If $V_1=\emptyset$, then it
follows from $n\geq 6$ that $V_2\neq\emptyset$. If
$V_1\neq\emptyset$, then by Lemma~\ref{lem3.2}, we also have
$V_2\neq\emptyset$.)

\begin{claim}\label{claim4.1}
$G[V_2]$ is a matching.
\end{claim}

\pf If there exist three vertices $u_1$, $u_2$ and $u_3$ in $V_2$
such that $u_1u_2,u_2u_3\in E(G)$, then $G[\{u_1,u_2,u_3,x_1,x_2\}]$
contains a copy of $W_4$, a contradiction. Thus, we conclude that
every component of $G[V_2]$ contains at most two vertices. Let $u$
be an arbitrary vertex in $V_2$. Because $G$ is $W_4$-saturated,
there exists a copy of $W_4$ in $G+ux$, say $H$. Since
$N_G(x)=\{x_1,x_2\}$, we know that $u,x,x_1,x_2\in V(H)$ and the
center of $H$ is $u$, $x_1$ or $x_2$. Let $u'$ be the remaining
vertex of $V(H)\setminus\{u,x,x_1,x_2\}$. It is easy to verify that
no matter the center of $H$ is $u$, $x_1$ or $x_2$, we always have
$uu',u'x_1,u'x_2\in E(G)$. This implies that every component of
$G[V_2]$ is a $K_2$, i.e. $G[V_2]$ is a matching.   \qed

\bigskip

By (1) and Claim~\ref{claim4.1}, we see that $f(u)=2.5$ for any
$u\in V_2$.

\begin{claim}\label{claim4.2}
$V_1\neq\emptyset$.
\end{claim}

\pf Suppose to the contrary that $V_1=\emptyset$. Then by
Claim~\ref{claim4.1}, we derive that $d(u)=3$ for any $u\in V_2$. If
$x_1x_2\in E(G)$, then both $x_1$ and $x_2$ are universal vertices
of $G$, contradicting the assumption that $G$ contains no universal
vertex. Hence, $x_1x_2\notin E(G)$. Let $H$ be a copy of $W_4$ in
$G+x_1x_2$. Since $x_1$ and $x_2$ are the only two vertices of
$G+x_1x_2$ with degree at least $4$, we may assume by symmetry that
$x_1$ is the center of $H$. Let $x_2u_1u_2u_3x_2$ be the rim of $H$.
Since $d(x)=2$ and $V_1=\emptyset$, we notice that $u_1,u_2,u_3\in
V_2$. But this contradicts Claim~\ref{claim4.1} (because
$u_1u_2,u_2u_3\in E(G)$). \qed

\bigskip

By Claim~\ref{claim4.2} and Lemma~\ref{lem3.2}, we have $x_1x_2\in
E(G)$. Let $V_1^*$ be the set of vertices in $V_1$ with degree $2$.
Then by (1) and Lemma~\ref{lem3.2}, we deduce that $f(v)=2$ for any
$v\in V_1^*$ and $f(v)\geq 2.5$ for any $v\in V_1\setminus V_1^*$.

\begin{claim}\label{claim4.3}
$|V_1^*|\leq 2$.
\end{claim}

\pf Suppose not, and let $v_1$, $v_2$ and $v_3$ be three vertices in
$V_1^*$. Without loss of generality, we may assume that
$v_1x_1,v_2x_1\in E(G)$. For each $i\in [3]$, let $v_i'$ be a shadow
of $v_i$ in $V_2$ (by Lemma~\ref{lem3.2}). Since $d(v_i)=2$ and
$N(v_i)\cap N(x)\neq\emptyset$ for each $i\in [3]$, we have
$v_1v_2,v_1v_3,v_2v_3\notin E(G)$. Let $H$ be a copy of $W_4$ in
$G+v_1v_2$. Then, it is easy to observe that the center of $H$ must
be $x_1$ and the rim of $H$ must be $v_1v_2v_2'v_1'v_1$. This
implies that $v_1'\neq v_2'$ and $v_1'v_2'\in E(G)$. If $v_3x_1\in
E(G)$, then by considering the copies of $W_4$ in $G+v_1v_3$ and in
$G+v_2v_3$, respectively, we can also conclude that $v_3'\notin
\{v_1',v_2'\}$ and $v_1'v_3',v_2'v_3'\in E(G)$, which contradicts
Claim~\ref{claim4.1}. Therefore, we have $v_3x_2\in E(G)$. Then by
Lemma~\ref{lem3.1}(ii), we know that $v_3'=v_1'$; otherwise, $v_3$
and $v_1$ have no common neighbor. But now, we see that $v_3$ and
$v_2$ have no common neighbor (since $v_1'\neq v_2'$), contradicting
Lemma~\ref{lem3.1}(ii).   \qed

\bigskip

By Lemma~\ref{lem3.4} and Claim~\ref{claim4.3}, we have
\begin{align*}
e(G)\geq 3+2\cdot 2+2.5(n-5)=\frac{5n-11}{2},
\end{align*}
i.e. $e(G)\geq \lfloor\frac{5n-10}{2}\rfloor$.

In the following, we characterize the extremal graphs. Suppose
$e(G)=\lfloor\frac{5n-10}{2}\rfloor$. If $V_1^*=\emptyset$, then it
follows from Lemma~\ref{lem3.4} that
\begin{align*}
e(G)\geq 3+2.5(n-3)=\frac{5n-9}{2}>\lfloor\frac{5n-10}{2}\rfloor,
\end{align*}
a contradiction. Thus by Claim~\ref{claim4.3}, we derive that $1\leq
|V_1^*|\leq 2$.

Let $y$ be a vertex in $V_1^*$. By symmetry between $x_1$ and $x_2$,
we may assume that $N(y)=\{x_1,u_1\}$, where $u_1$ is the shadow of
$y$ in $V_2$ (by Lemma~\ref{lem3.2}). By Claim~\ref{claim4.1}, let
$u_2$ be the unique neighbor of $u_1$ in $V_2$. Define
$U_2:=V_2\setminus\{u_1,u_2\}$. Then by Claim~\ref{claim4.1} and the
definition of $V_2$, we deduce that if $U_2\neq\emptyset$ then
$G[U_2]$ is a matching and every vertex in $U_2$ is adjacent to both
$x_1$ and $x_2$. For each edge $vv'\in E(G[V_1])$ with
$d(v)=d(v')=3$, we say that $vv'$ is of \emph{Type 1} if $N(v)\cap
N(x)=N(v')\cap N(x)$ and of \emph{Type 2} if $N(v)\cap N(x)\neq
N(v')\cap N(x)$.

\begin{claim}\label{claim4.4}
If $v_1v_2$ is a Type $2$ edge in $G[V_1]$ such that
$v_1x_1,v_2x_2\in E(G)$, then $v_1u_2\in E(G)$.
\end{claim}

\pf First, suppose $N(v_1)\cap\{u_1,u_2\}=\emptyset$. Since
$d(v_1)=3$, we may assume that $N(v_1)=\{x_1,v_2,u_3\}$, where $u_3$
is the shadow of $v_1$ in $V_2$ (by Lemma~\ref{lem3.2}). Then by
Claim~\ref{claim4.1}, we have $u_3u_1,u_3u_2\notin E(G)$. Let $H$ be
a copy of $W_4$ in $G+v_1y$. It is easy to see that the center of
$H$ must be $x_1$ (since $x_1$ is the unique common neighbor of
$v_1$ and $y$ in $G$) and the rim of $H$ must be $v_1yu_1v_2v_1$.
But this implies that $v_2x_1\in E(G)$, which contradicts the
assumption that $v_2x_1\notin E(G)$.

Hence, $N(v_1)\cap\{u_1,u_2\}\neq\emptyset$. If $v_1u_1\in E(G)$,
then we know that $N(v_1)=\{x_1,v_2,u_1\}$ (since $d(v_1)=3$) and it
is easy to verify that $G+v_1y$ contains no copy of $W_4$ (since
$v_2x_1\notin E(G)$), again a contradiction. Therefore, we have $v_1u_2\in
E(G)$. \qed

\begin{claim}\label{claim4.5}
$G[V_1]$ contains at most one Type $2$ edge.
\end{claim}

\pf Suppose not, and let $v_1v_2,v_3v_4$ be two Type 2 edges in
$G[V_1]$. Without loss of generality, we may assume that
$v_1x_1,v_2x_2,v_3x_1,v_4x_2\in E(G)$. Then by Claim~\ref{claim4.4},
we see that $v_1u_2,v_3u_2\in E(G)$ and thus
$N(v_3)=\{x_1,v_4,u_2\}$. This means that $v_2u_2\in E(G)$;
otherwise, $v_2$ and $v_3$ have no common neighbor, contradicting
Lemma~\ref{lem3.1}(ii). But now, $G[\{u_2,v_1,v_2,x_1,x_2\}]$
contains a copy of $W_4$, a contradiction. \qed

\bigskip

We now consider two cases according to the value of $|V_1^*|$.

\medskip

{\bf Case 1.} $|V_1^*|=1$.

\medskip

In this case, $y$ is the unique vertex in $V_1^*$ and
$E(G[V_1])=E(G[V_1\setminus V_1^*])$. It is clear that $V_1\setminus
V_1^*\neq\emptyset$; otherwise, $x_1$ would be a universal vertex of
$G$, contradicting the assumption that $G$ contains no universal
vertex. If there exists some vertex $v\in V_1\setminus V_1^*$ such
that $f(v)\geq 3$, then it follows from Lemma~\ref{lem3.4} that
\begin{align*}
e(G)\geq
3+2+3+2.5(n-5)=\frac{5n-9}{2}>\lfloor\frac{5n-10}{2}\rfloor,
\end{align*}
a contradiction. Thus, we have $f(v)=2.5$ for any $v\in
V_1\setminus V_1^*$. This implies that $d(v)=3$ for any $v\in
V_1\setminus V_1^*$ and $G[V_1\setminus V_1^*]$ is a matching (by
(1) and Lemma~\ref{lem3.2}). Note that $G[V_1\setminus V_1^*]$
contains at most one Type $2$ edge by Claim~\ref{claim4.5}.

\medskip

{\bf Subcase 1.1.} $G[V_1\setminus V_1^*]$ contains one Type $2$
edge, say $v_1v_2$.

\medskip

Without loss of generality, suppose $v_1x_1,v_2x_2\in E(G)$. Then by
Claim~\ref{claim4.4}, we derive that $v_1u_2\in E(G)$. Moreover, we
have $v_2u_1\in E(G)$; otherwise, $v_2$ and $y$ have no common
neighbor, contradicting Lemma~\ref{lem3.1}(ii). This shows that
$N(v_1)=\{x_1,v_2,u_2\}$ and $N(v_2)=\{x_2,v_1,u_1\}$. Define
$U_1:=V_1\setminus\{y,v_1,v_2\}$. If $U_1=\emptyset$, then we can
deduce that $G\in\mathcal{B}_n^3$. So we may assume that
$U_1\neq\emptyset$. Then $G[U_1]$ is still a matching. If there
exists some vertex $v\in U_1$ such that $vx_2\in E(G)$, then by
Lemma~\ref{lem3.1}(ii), we conclude that $vu_1\in E(G)$; otherwise,
$v$ and $y$ have no common neighbor. But then, $v$ and $v_1$ have no
common neighbor, which contradicts Lemma~\ref{lem3.1}(ii). Hence,
every vertex in $U_1$ is adjacent to $x_1$. On the other hand, we
notice that every vertex in $U_1$ is also adjacent to $u_1$;
otherwise, $v$ and $v_2$ have no common neighbor for some vertex
$v\in U_1$, contradicting Lemma~\ref{lem3.1}(ii). Then, it is
straightforward to check that $G\in\mathcal{B}_n^3$.

\medskip

{\bf Subcase 1.2.} $G[V_1\setminus V_1^*]$ contains no Type $2$
edge.

\medskip

Then every edge in $G[V_1\setminus V_1^*]$ is of Type $1$. For each
$i\in [2]$, let $U_1^i$ be the set of vertices in $V_1\setminus
V_1^*$ that are adjacent to $x_i$. Since $G$ contains no universal
vertex, we know that $U_1^2\neq\emptyset$; otherwise, $x_1$ would be
a universal vertex of $G$, a contradiction. Moreover, it follows
from $G[V_1\setminus V_1^*]$ is a matching that both $G[U_1^1]$ (if
$U_1^1\neq\emptyset$) and $G[U_1^2]$ are matchings. It is clear that
every vertex in $U_1^2$ is adjacent to $u_1$; otherwise, $v$ and $y$
have no common neighbor for some $v\in U_1^2$, which contradicts
Lemma~\ref{lem3.1}(ii). If $U_1^1=\emptyset$, then we see that
$G\in\mathcal{B}_n^2$. Therefore, we may assume that
$U_1^1\neq\emptyset$. If there exists some vertex $v'\in U_1^1$ such
that $v'u_1\notin E(G)$, then $v'$ and $v$ have no common neighbor
for any $v\in U_1^2$, contradicting Lemma~\ref{lem3.1}(ii). Thus,
every vertex in $U_1^1$ is adjacent to $u_1$. This also implies that
$G\in\mathcal{B}_n^2$.

\medskip

{\bf Case 2.} $|V_1^*|=2$.

\medskip

Let $z$ be the other vertex in $V_1^*$ except $y$. Then $N(y)\neq
N(z)$; otherwise, one can easily check that $G+yz$ contains no copy
of $W_4$, a contradiction.

Suppose $zx_1\in E(G)$. Let $H$ be a copy of $W_4$ in $G+yz$. Then
by Lemma~\ref{lem3.2} and Claim~\ref{claim4.1}, we can derive that
the center of $H$ must be $x_1$ and the rim of $H$ must be
$yzu_2u_1y$. This means that $N(z)=\{x_1,u_2\}$. Since $G$ contains
no universal vertex, we deduce that there must exist a vertex $v\in
V_1\setminus V_1^*$ such that $vx_2\in E(G)$; otherwise, $x_1$ would
be a universal vertex of $G$, a contradiction. Then by
Lemma~\ref{lem3.1}(ii), we conclude that $vu_1,vu_2\in E(G)$;
otherwise, either $v$ and $y$ (if $vu_1\notin E(G)$) or $v$ and $z$
(if $vu_2\notin E(G)$) have no common neighbor. But then,
$G[\{v,u_1,u_2,x_1,x_2\}]$ contains a copy of $W_4$, giving a
contradiction.

Therefore, we have $zx_2\in E(G)$. Then by Lemma~\ref{lem3.1}(ii),
we know that $N(z)=\{x_2,u_1\}$; otherwise, $z$ and $y$ have no
common neighbor. If $V_1\setminus V_1^*=\emptyset$, then we see that
$G\in\mathcal{B}_n^1$. Hence, we may assume that $V_1\setminus
V_1^*\neq\emptyset$. Then, it is easy to observe that every vertex
in $V_1\setminus V_1^*$ is adjacent to $u_1$; otherwise, either $v$
and $z$ (if $vx_1\in E(G)$) or $v$ and $y$ (if $vx_2\in E(G)$) have
no common neighbor for some vertex $v\in V_1\setminus V_1^*$,
contradicting Lemma~\ref{lem3.1}(ii).

If there exists some vertex $v\in V_1\setminus V_1^*$ such that
$f(v)\geq 3.5$ or two vertices $v,v'\in V_1\setminus V_1^*$ such
that $f(v)=f(v')=3$, then by Lemma~\ref{lem3.4}, we have
\begin{align*}
e(G)\geq 3+2\cdot
2+3.5+2.5(n-6)=\frac{5n-9}{2}>\lfloor\frac{5n-10}{2}\rfloor
\end{align*}
or
\begin{align*}
e(G)\geq 3+2\cdot 2+3\cdot
2+2.5(n-7)=\frac{5n-9}{2}>\lfloor\frac{5n-10}{2}\rfloor,
\end{align*}
a contradiction. Thus, we derive that $2.5\leq f(v)\leq
3$ for any $v\in V_1\setminus V_1^*$ and there is at most one vertex
$v'\in V_1\setminus V_1^*$ such that $f(v')=3$.

\medskip

{\bf Subcase 2.1.} There exists a vertex $v'\in V_1\setminus
V_1^*$ such that $f(v')=3$.

\medskip

Without loss of generality, we may assume that $v'x_1\in E(G)$.
Since every vertex in $V_1\setminus V_1^*$ is adjacent to $u_1$, we
have $v'u_1\in E(G)$. If there exists another vertex $u'\in
V_2\setminus\{u_1\}$ such that $v'u'\in E(G)$, then
$G[\{v',u',u_1,x_1,x_2\}]$ contains a copy of $W_4$, a
contradiction. Hence, $u_1$ is the unique neighbor of $v'$ in $V_2$.
Since $f(v')=3$ and by (1), we deduce that there must exist two
vertices $v_1,v_2\in V_1\setminus V_1^*$ such that $v'v_1,v'v_2\in
E(G)$. Note that $v_1u_1,v_2u_1\in E(G)$ (since every vertex in
$V_1\setminus V_1^*$ is adjacent to $u_1$). If there exists some
$i\in [2]$ such that $v_1x_i,v_2x_i\in E(G)$, then
$G[\{v',v_1,v_2,u_1,x_i\}]$ contains a copy of $W_4$, giving a
contradiction. Therefore, we may assume by symmetry that
$v_1x_1,v_2x_2\in E(G)$. But then, $G[\{v',u_1,v_2,x_1,x_2\}]$
contains a copy of $W_4$, a contradiction.

\medskip

{\bf Subcase 2.2.} There is no vertex $v'\in V_1\setminus V_1^*$
such that $f(v')=3$.

\medskip

Then $f(v)=2.5$ for any $v\in V_1\setminus V_1^*$. By (1) and
Lemma~\ref{lem3.2}, we conclude that $d(v)=3$ for any $v\in
V_1\setminus V_1^*$ and $G[V_1\setminus V_1^*]$ is a matching.

Suppose $G[V_1\setminus V_1^*]$ contains a Type $2$ edge, say
$v_1v_2$. Without loss of generality, we may assume that
$v_1x_1,v_2x_2\in E(G)$. Since every vertex in $V_1\setminus V_1^*$
is adjacent to $u_1$, we know that $v_1u_1\in E(G)$. On the other
hand, it follows from Claim~\ref{claim4.4} that $v_1u_2\in E(G)$.
But this implies that $f(v_1)\geq 3$ (by (1)), a contradiction.

Thus, we see that $G[V_1\setminus V_1^*]$ contains no Type $2$ edge.
For each $i\in [2]$, let $U_1^i$ be the set of vertices in
$V_1\setminus V_1^*$ that are adjacent to $x_i$. Since
$G[V_1\setminus V_1^*]$ is a matching and every vertex in
$V_1\setminus V_1^*$ is adjacent to $u_1$, we can derive that for
each $i\in [2]$, if $U_1^i\neq\emptyset$ then $G[U_1^i]$ is also a
matching and every vertex in $U_1^i$ is adjacent to $u_1$. Then, it
is straightforward to verify that $G\in\mathcal{B}_n^1$.

\medskip

This completes the proof of the first part. \qed

\subsection{$\delta(G)=3$ and $e(N[x])=3$}

\medskip

In this part, $V_1=\emptyset$ (by Lemma~\ref{lem3.2}).

\begin{claim}\label{claim4.6}
For any $u\in V_2\cup V_3$, there exists a vertex $w\in V_2\cup
V_3$ such that $uw\in E(G)$ and $|N(u)\cap N(w)\cap N(x)|\geq 2$.
\end{claim}

\pf Let $H$ be a copy of $W_4$ in $G+ux$. Since $e(N[x])=3$, we
deduce that the center of $H$ is $u$. Let $xx_iwx_jx$ be the rim of
$H$ for some $i,j\in [3]$. Then $uw\in E(G)$ and
$\{x_i,x_j\}\subseteq N(u)\cap N(w)$. This implies that $w\in
V_2\cup V_3$ and $|N(u)\cap N(w)\cap N(x)|\geq 2$.      \qed

\begin{claim}\label{claim4.7}
For any $u\in V_2$ with $d(u)=3$, we have $N(u)\cap
V_3\neq\emptyset$.
\end{claim}

\pf Suppose to the contrary that $N(u)\cap V_3=\emptyset$ for some
vertex $u\in V_2$ with $d(u)=3$. Without loss of generality, we may
assume that $N(u)=\{x_1,x_2,w\}$, where $w$ is the unique neighbor
of $u$ in $V_2$. Then by Claim~\ref{claim4.6}, we conclude that
$wx_1,wx_2\in E(G)$. But then, $u$ and $x_3$ have no common
neighbor, contradicting Lemma~\ref{lem3.1}(ii).     \qed

\bigskip

It follows from Claim~\ref{claim4.6} that every component of
$G[V_2\cup V_3]$ contains at least two vertices, and thus $d(w)\geq
4$ for any $w\in V_3$. Then by (2), we know that $g(u)\geq 2.5$ for
any $u\in V_2$ and $g(w)\geq 3.5$ for any $w\in V_3$. If $|V_3|\geq
3$, then by Lemma~\ref{lem3.5}, we have
\begin{align*}
e(G)\geq 3+3.5\cdot
3+2.5(n-7)=\frac{5n-8}{2}>\lfloor\frac{5n-10}{2}\rfloor.
\end{align*}
Hence, we may assume that $|V_3|\leq 2$.

We now consider three cases according to the value of $|V_3|$.

\medskip

{\bf Case 1.} $|V_3|=0$.

\medskip

Since $n\geq 6$, we see that $V_2\neq\emptyset$. Then by
Claim~\ref{claim4.7}, we have $d(u)\geq 4$ for any $u\in V_2$. This
shows that $g(u)\geq 3$ for any $u\in V_2$ (by (2)) and $|V_2|\geq
3$. If $n\geq 9$, then we derive that
\begin{align*}
e(G)\geq
3+3(n-4)=3n-9=\frac{5n-9}{2}+\frac{n-9}{2}\geq\frac{5n-9}{2}>\lfloor\frac{5n-10}{2}\rfloor
\end{align*}
by Lemma~\ref{lem3.5}. Therefore, we may assume that $n\leq 8$ and
thus $|V_2|\leq 4$. Then $3\leq |V_2|\leq 4$.

First, suppose $|V_2|=3$ (i.e. $n=7$). Let $V_2=\{u_1,u_2,u_3\}$.
Since $d(u_i)\geq 4$ for each $i\in [3]$, we deduce that
$G[V_2]\cong C_3$. Without loss of generality, suppose
$u_1x_1,u_1x_2\in E(G)$. Then by Claim~\ref{claim4.6}, we may assume
by symmetry that $u_2x_1,u_2x_2\in E(G)$. But this implies that
$d(x_3)\leq 2$, which contradicts the assumption that $\delta(G)=3$.

Next, suppose $|V_2|=4$ (i.e. $n=8$). Let $V_2=\{u_1,u_2,u_3,u_4\}$.
If there exists some $i\in [4]$ such that $d(u_i)\geq 5$, then
$g(u_i)\geq 3.5$ (by (2)) and it follows from Lemma~\ref{lem3.5}
that
\begin{align*}
e(G)\geq 3+3.5+3\cdot 3=15.5>15=\lfloor\frac{5n-10}{2}\rfloor.
\end{align*}
Thus, we may assume that $d(u_i)=4$ for each $i\in [4]$, and hence
$G[V_2]\cong C_4$. Let $G[V_2]=u_1u_2u_3u_4u_1$ and suppose without
loss of generality that $u_1x_1,u_1x_2\in E(G)$. Then by
Claim~\ref{claim4.6}, we may assume by symmetry that
$u_2x_1,u_2x_2\in E(G)$. This means that $u_3x_3,u_4x_3\in E(G)$;
otherwise, either $u_2$ and $x_3$ (if $u_3x_3\notin E(G)$) or $u_1$
and $x_3$ (if $u_4x_3\notin E(G)$) have no common neighbor,
contradicting Lemma~\ref{lem3.1}(ii). By symmetry between $x_1$ and
$x_2$, we may further assume that $u_3x_1\in E(G)$. Then by
Claim~\ref{claim4.6}, we can conclude that $u_4x_1\in E(G)$. But
now, it is easy to check that $G+x_2x_3$ contains no copy of $W_4$
(since $x$ is the unique common neighbor of $x_2$ and $x_3$ in $G$
and $d(x)=3$), a contradiction.

\medskip

{\bf Case 2.} $|V_3|=1$.

\medskip

In this case, we also have $V_2\neq\emptyset$ (since $n\geq 6$). Let
$V_3=\{w\}$.

First, suppose there exists a vertex $u_1\in V_2$ such that
$u_1w\notin E(G)$. Then by Claim~\ref{claim4.7}, we notice that
$d(u_1)\geq 4$. Let $u_2$ and $u_3$ be two neighbors of $u_1$ in
$V_2$. Then $d(u_2)\geq 4$ and $d(u_3)\geq 4$. (For each $i\in
\{2,3\}$, if $u_iw\in E(G)$ then it is clear that $d(u_i)\geq 4$,
and if $u_iw\notin E(G)$ then we also have $d(u_i)\geq 4$ by
Claim~\ref{claim4.7}.) By (2), we know that $g(u_i)\geq 3$ for each
$i\in [3]$. Then, it follows from Lemma~\ref{lem3.5} that
\begin{align*}
e(G)\geq 3+3.5+3\cdot
3+2.5(n-8)=\frac{5n-9}{2}>\lfloor\frac{5n-10}{2}\rfloor.
\end{align*}

Next, suppose every vertex in $V_2$ is adjacent to $w$. Let $u_1$ be
a neighbor of $w$ in $V_2$, and assume without loss of generality
that $u_1x_1,u_1x_2\in E(G)$. Let $H$ be a copy of $W_4$ in
$G+u_1x_3$. One can easily check that no matter the center of $H$ is
$u_1$, $x_3$ or some common neighbor of $u_1$ and $x_3$, there must
exist a vertex $u_2\in V(H)\cap V_2$ such that $u_2x_3,u_2w\in
E(G)$. This implies that $d(w)\geq 5$, and thus $g(w)\geq 4$ (by
(2)). By symmetry between $x_1$ and $x_2$, we may assume that
$u_2x_1\in E(G)$. If $d(u_1)\geq 4$ and $d(u_2)\geq 4$, then
$g(u_i)\geq 3$ for each $i\in [2]$ (by (2)) and
\begin{align*}
e(G)\geq 3+4+3\cdot
2+2.5(n-7)=\frac{5n-9}{2}>\lfloor\frac{5n-10}{2}\rfloor
\end{align*}
by Lemma~\ref{lem3.5}. Hence, we may assume by symmetry that
$d(u_1)=3$ (and thus $u_1u_2\notin E(G)$). Then we see that
$d(u_2)\geq 4$; otherwise, it is easy to verify that $G+u_1u_2$
contains no copy of $W_4$ (since $x_1x_2,x_1x_3,x_2x_3\notin E(G)$),
a contradiction. This shows that $|V_2|\geq 3$. Since every vertex
in $V_2$ is adjacent to $w$, we derive that $d(w)\geq 6$. Then by
(2), we have $g(u_2)\geq 3$ and $g(w)\geq 4.5$. Now, it follows from
Lemma~\ref{lem3.5} that
\begin{align*}
e(G)\geq
3+4.5+3+2.5(n-6)=\frac{5n-9}{2}>\lfloor\frac{5n-10}{2}\rfloor.
\end{align*}

\medskip

{\bf Case 3.} $|V_3|=2$.

\medskip

Let $V_3=\{w_1,w_2\}$. Suppose $V_2=\emptyset$. Then by
Claim~\ref{claim4.6}, we deduce that $w_1w_2\in E(G)$ and hence
$d(w_1)=d(w_2)=4$. But then, since $x_1x_3,x_2x_3\notin E(G)$, it is
straightforward to check that $G+x_1x_2$ contains no copy of $W_4$,
a contradiction. Therefore, we have $V_2\neq\emptyset$.

Recall that $d(w_i)\geq 4$ and $g(w_i)\geq 3.5$ for each $i\in [2]$
(by Claim~\ref{claim4.6} and (2)). If there exists some vertex $u\in
V_2$ such that $d(u)\geq 4$, then $g(u)\geq 3$ (by (2)) and
\begin{align*}
e(G)\geq 3+3.5\cdot
2+3+2.5(n-7)=\frac{5n-9}{2}>\lfloor\frac{5n-10}{2}\rfloor
\end{align*}
by Lemma~\ref{lem3.5}. Similarly, if there exists some
$i\in [2]$ such that $d(w_i)\geq 5$, then $g(w_i)\geq 4$ (by (2))
and it follows from Lemma~\ref{lem3.5} that
\begin{align*}
e(G)\geq
3+4+3.5+2.5(n-6)=\frac{5n-9}{2}>\lfloor\frac{5n-10}{2}\rfloor.
\end{align*}
Thus, we may assume that $d(u)=3$ for any $u\in V_2$ and $d(w_i)=4$
for each $i\in [2]$. This implies that $G[V_2\cup V_3]$ is a
matching. Moreover, we observe that every vertex in $V_2$ is
adjacent to a vertex in $V_3$ (by Claim~\ref{claim4.7}), and hence
$|V_2|=2$. Let $V_2=\{u_1,u_2\}$ such that $u_1w_1,u_2w_2\in E(G)$.
Without loss of generality, we may assume that
$N(u_1)=\{x_1,x_2,w_1\}$ and $N(u_2)=\{x_1,x_j,w_2\}$ for some
$j\in\{2,3\}$. Now, one can easily see that $G+u_1x_3$ contains no
copy of $W_4$ (since $x_1x_2,x_1x_3,x_2x_3\notin E(G)$), a
contradiction.

\medskip

In conclusion, we show that $e(G)>\lfloor\frac{5n-10}{2}\rfloor$ in
all cases and there is no extremal graph in this part.  \qed

\subsection{$\delta(G)=3$ and $e(N[x])=4$}

\medskip

In this part, suppose without loss of generality that $x_1x_2\in
E(G)$ and $x_1x_3,x_2x_3\notin E(G)$.

\begin{claim}\label{claim4.8}
The following statements hold:
\begin{itemize}
\item [$($i$)$] $vx_3\notin E(G)$ for any $v\in V_1$;
\item [$($ii$)$] $N(u)\cap (V_2\cup V_3)\neq\emptyset$ for any $u\in V_2$;
\item [$($iii$)$] if $V_3=\emptyset$, then $ux_3\in E(G)$ for any $u\in V_2$
with $d(u)=3$.
\end{itemize}
\end{claim}

\pf Let $v$ be a vertex in $V_1$ such that $vx_3\in E(G)$. Then by
Lemma~\ref{lem3.2}, we conclude that $x_jx_3\in E(G)$ for some $j\in
[2]$, which contradicts the assumption that $x_1x_3,x_2x_3\notin
E(G)$. This proves (i).

Let $u$ be a vertex in $V_2$ such that $N(u)\cap (V_2\cup
V_3)=\emptyset$. Then by Lemma~\ref{lem3.3}, we know that
$e(N[x])\geq 5$, contradicting the assumption that $e(N[x])=4$. So
we have (ii).

Finally, we prove (iii). Suppose not, and let $u$ be a vertex in
$V_2$ such that $d(u)=3$ and $ux_3\notin E(G)$. Since $V_3=\emptyset$ and by
Claim~\ref{claim4.8}(ii), we may assume that $N(u)=\{x_1,x_2,u'\}$,
where $u'$ is the unique neighbor of $u$ in $V_2$. Then by
Lemma~\ref{lem3.1}(ii), we see that $u'x_3\in E(G)$; otherwise, $u$
and $x_3$ have no common neighbor. Let $H$ be a copy of $W_4$ in
$G+ux$. Since $x_1x_3,x_2x_3\notin E(G)$, we notice that $x$ is not
the center of $H$. This means that the center of $H$ must be one
vertex in $\{u,x_1,x_2\}$. It is easy to check that in all cases, we
always have $V(H)=\{u,u',x,x_1,x_2\}$ and $u'x_1,u'x_2\in E(G)$. But
this implies that $u'\in V_3$, contradicting the assumption that
$V_3=\emptyset$.   \qed

\bigskip

It follows from Lemma~\ref{lem3.2} and Claim~\ref{claim4.8}(ii) that
every vertex in $V_1\cup V_2$ has at least one neighbor in $V_2\cup
V_3$. Then by (1), we derive that $f(v)\geq 2.5$ for any $v\in
V_1\cup V_2$ and $f(w)\geq 3$ for any $w\in V_3$.

\begin{claim}\label{claim4.9}
If $V_3=\emptyset$, then there exists a vertex $v\in V_1\cup V_2$
such that $f(v)\geq 3$.
\end{claim}

\pf Let $H$ be a copy of $W_4$ in $G+x_2x_3$ and let $z$ be the
center of $H$. Then $z$ is $x_2$, $x_3$ or some common neighbor of
$x_2$ and $x_3$.

First, suppose $z\in\{x_2,x_3\}$. We only deal with the case that
$z=x_2$ here, while the case that $z=x_3$ can be treated in a
similar way. Let $x_3y_1y_2y_3x_3$ be the rim of $H$. If
$\{y_1,y_2,y_3\}\cap\{x,x_1\}\neq\emptyset$, then it is easy to
observe that $y_2=x_1$ and $y_i\in V_3$ for some $i\in\{1,3\}$,
contradicting the assumption that $V_3=\emptyset$. Hence, we deduce
that $y_1,y_3\in V_2$ (since $y_1x_2,y_1x_3,y_3x_2,y_3x_3\in E(G)$)
and $y_2\in V_1\cup V_2$. Since $y_2y_1,y_2y_3\in E(G)$ and by (1),
we have $f(y_2)\geq 3$ (no matter $y_2\in V_1$ or $y_2\in V_2$), as
desired.

Next, suppose $z$ is some common neighbor of $x_2$ and $x_3$. Let
$x_2x_3y_1y_2x_2$ be the rim of $H$. Since $d(x)=3$ and
$x_1x_3\notin E(G)$, we conclude that $z\notin\{x,x_1\}$. This shows
that $z\in V_2$ (since $zx_2,zx_3\in E(G)$ and $V_3=\emptyset$), and
thus $y_1,y_2\in V_1\cup V_2$. Moreover, because $y_1x_3\in E(G)$,
it follows from Claim~\ref{claim4.8}(i) that $y_1\in V_2$. Note that
$y_2y_1,y_2z\in E(G)$. Then by (1), we know that $f(y_2)\geq 3$ (no
matter $y_2\in V_1$ or $y_2\in V_2$), as required.     \qed

\bigskip

By Claim~\ref{claim4.9} and (1), we see that no matter whether
$V_3=\emptyset$ or not, there always exists a vertex $v\in V_x$ such
that $f(v)\geq 3$. Then by Lemma~\ref{lem3.4}, we have
\begin{align*}
e(G)\geq 4+3+2.5(n-5)=\frac{5n-11}{2},
\end{align*}
i.e. $e(G)\geq \lfloor\frac{5n-10}{2}\rfloor$.

In the following, we characterize the extremal graphs. Suppose
$e(G)=\lfloor\frac{5n-10}{2}\rfloor$. If there exists some vertex
$v\in V_x$ such that $f(v)\geq 4$, then by Lemma~\ref{lem3.4}, we
derive that
\begin{align*}
e(G)\geq
4+4+2.5(n-5)=\frac{5n-9}{2}>\lfloor\frac{5n-10}{2}\rfloor,
\end{align*}
a contradiction. Therefore, we have $2.5\leq f(v)\leq 3.5$ for any
$v\in V_x$. Let $V_x^*:=\{v\in V_x: 3\leq f(v)\leq 3.5\}$. Then
$V_3\subseteq V_x^*$.

\begin{claim}\label{claim4.10}
$|V_x^*|\leq 2$. Moreover, if $|V_x^*|=2$, then $f(v)=3$ for any
$v\in V_x^*$.
\end{claim}

\pf If $|V_x^*|\geq 3$, then it follows from Lemma~\ref{lem3.4} that
\begin{align*}
e(G)\geq 4+3\cdot
3+2.5(n-7)=\frac{5n-9}{2}>\lfloor\frac{5n-10}{2}\rfloor,
\end{align*}
a contradiction. Thus, $|V_x^*|\leq 2$.

If $|V_x^*|=2$ and there exists some vertex $v\in V_x^*$ such that
$f(v)=3.5$, then by Lemma~\ref{lem3.4}, we deduce that
\begin{align*}
e(G)\geq
4+3.5+3+2.5(n-6)=\frac{5n-9}{2}>\lfloor\frac{5n-10}{2}\rfloor,
\end{align*}
again a contradiction. \qed

\begin{claim}\label{claim4.11}
If $V_1\neq\emptyset$, then $V_x^*\cap V_1\neq\emptyset$.
\end{claim}

\pf Suppose to the contrary that $V_x^*\cap V_1=\emptyset$. Then
$f(v)=2.5$ for any $v\in V_1$. By (1) and Lemma~\ref{lem3.2}, we
conclude that $d(v)=3$ for any $v\in V_1$ and $G[V_1]$ is a matching
(i.e. every vertex in $V_1$ has exactly one neighbor in $V_1$ and
exactly one neighbor in $V_2\cup V_3$). Let $v_1v_2$ be an edge in
$G[V_1]$. Then by Claim~\ref{claim4.8}(i), we have
$v_1x_3,v_2x_3\notin E(G)$. Since $d(v_1)=3$, we may assume by
symmetry that $N(v_1)=\{x_1,v_2,w\}$, where $w$ is the shadow of
$v_1$ in $V_2\cup V_3$ and $wx_1,wx_2\in E(G)$ (by
Lemma~\ref{lem3.2}). Then by Lemma~\ref{lem3.1}(ii), we know that
$wx_3\in E(G)$; otherwise, $v_1$ and $x_3$ have no common neighbor.
This implies that $w\in V_3$.

Let $H_1$ be a copy of $W_4$ in $G+v_1x_3$. Then the center of $H_1$
must be $w$ (since $w$ is the unique common neighbor of $v_1$ and
$x_3$ in $G$). Let $v_1x_3yzv_1$ be the rim of $H_1$. Since $v_1z\in
E(G)$, we have $z\in\{x_1,v_2\}$. If $z=v_2$, then we see that
$v_2w\in E(G)$ and thus $y=x_i$ for some $i\in [2]$ (since
$d(v_2)=3$ and $v_2x_3\notin E(G)$), which contradicts the
assumption that $x_1x_3,x_2x_3\notin E(G)$. Hence, $z=x_1$. It is
clear that $y\notin\{x,x_2\}$. Note that $w\in V_3$ and $wy\in
E(G)$. If $y\in V_3$, then $f(w)=f(y)=3.5$ (by (1)), contradicting
Claim~\ref{claim4.10}. This shows that $y\in V_2$ (since
$yx_1,yx_3\in E(G)$). Then by Claim~\ref{claim4.10}, we can derive
that $V_x^*=\{w,y\}$ and $f(w)=f(y)=3$. Moreover, the following
statements hold:
\begin{itemize}
\item [(P1)] $V_3=\{w\}$ (since $V_x^*=\{w,y\}$ and $V_3\subseteq V_x^*$);
\item [(P2)] $vy\notin E(G)$ for any $v\in V_1$ (if $vy\in E(G)$
for some vertex $v\in V_1$, then $y$ is the shadow of $v$ in
$V_2\cup V_3$ and it follows from Lemma~\ref{lem3.2} that $x_1x_3\in
E(G)$, which contradicts the assumption that $x_1x_3\notin
E(G)$);
\item [(P3)] $uy\notin E(G)$ for any $u\in V_2$ (if $uy\in
E(G)$ for some vertex $u\in V_2$, then by (1), we have $f(y)= 3.5$,
contradicting the fact that $f(y)=3$);
\item [(P4)] $uw\notin E(G)$ for any $u\in V_2\setminus\{y\}$
(if $uw\in E(G)$ for some vertex $u\in V_2\setminus\{y\}$, then by
(1), we deduce that $f(u)\geq 3$, contradicting the fact that
$V_x^*=\{w,y\}$).
\end{itemize}

By (P1), (P2) and (P3), we conclude that $N(y)=\{x_1,x_3,w\}$. Let
$H_2$ be a copy of $W_4$ in $G+yx_2$. Since $x_1x_3,x_2x_3\notin
E(G)$, we observe that $y$ is not the center of $H_2$. This implies
that the center of $H_2$ is $w$, $x_1$ or $x_2$. But then, one can
easily check that in all cases, there must exist a vertex $u\in
V_2\setminus\{y\}$ such that $uw\in E(G)$, contradicting (P4). \qed

\begin{claim}\label{claim4.12}
$V_3=\emptyset$.
\end{claim}

\pf Suppose to the contrary that $V_3\neq\emptyset$. Since
$V_3\subseteq V_x^*$ and by Claim~\ref{claim4.10}, we have $1\leq
|V_3|\leq 2$.

First, suppose $|V_3|=1$. Let $V_3=\{w\}$. Because $G$ is
$W_4$-saturated, there exists a copy of $W_4$ in $G+wx$, say $H$.
Since $x_1x_3,x_2x_3\notin E(G)$, we know that neither $x$ nor $x_3$
is the center of $H$. Then the center of $H$ is $w$, $x_1$ or $x_2$.
It is easy to see that in all cases, there must exist a vertex $y\in
V_2$ such that $wy\in E(G)$. By (1) and Claim~\ref{claim4.10}, we
see that $V_x^*=\{w,y\}$ and $f(w)=f(y)=3$. This means that
$V_x^*\cap V_1=\emptyset$, and hence $V_1=\emptyset$ (by
Claim~\ref{claim4.11}). Moreover, by the same arguments as for (P3)
and (P4) in the proof of Claim~\ref{claim4.11}, we have $uy,uw\notin
E(G)$ for any $u\in V_2\setminus\{y\}$. Therefore, we derive that
$N(w)=\{x_1,x_2,x_3,y\}$ and $N(y)=\{x_i,x_j,w\}$ for some $i,j\in
[3]$. If $yx_3\notin E(G)$, then $N(y)=\{x_1,x_2,w\}$ and it is
straightforward to verify that $G+yx_3$ contains no copy of $W_4$
(since $x_1x_3,x_2x_3\notin E(G)$), a contradiction. Thus, $yx_3\in
E(G)$ and we may assume by symmetry that $N(y)=\{x_1,x_3,w\}$. But
then, since $x_1x_3,x_2x_3\notin E(G)$, it is easy to observe that
$G+yx_2$ contains no copy of $W_4$, giving a contradiction.

Next, suppose $|V_3|=2$. Let $V_3=\{w_1,w_2\}$. Then by (1) and
Claim~\ref{claim4.10}, we deduce that $V_x^*=\{w_1,w_2\}$ and
$f(w_1)=f(w_2)=3$. This implies that $w_1w_2\notin E(G)$ and
$V_1=\emptyset$ (by Claim~\ref{claim4.11}). Moreover, we notice that
$uw_1,uw_2\notin E(G)$ for any $u\in V_2$; otherwise, $f(u)\geq 3$
for some vertex $u\in V_2$ (by (1)), contradicting the fact that
$V_x^*=\{w_1,w_2\}$. Hence, we have $N(w_1)=N(w_2)=\{x_1,x_2,x_3\}$.
But now, one can easily see that $G+w_1w_2$ contains no copy of
$W_4$ (since $x_1x_3,x_2x_3\notin E(G)$), a contradiction.   \qed

\bigskip

By Claim~\ref{claim4.12}, we conclude that $V_2\neq\emptyset$. (If
$V_1=\emptyset$, then it follows from $n\geq 6$ that
$V_2\neq\emptyset$. If $V_1\neq\emptyset$, then by
Lemma~\ref{lem3.2}, we also have $V_2\neq\emptyset$.)

\begin{claim}\label{claim4.13}
If $V_1\neq\emptyset$, then $f(v)=3$ for any $v\in V_1$.
\end{claim}

\pf Let $v$ be an arbitrary vertex in $V_1$, and assume without loss
of generality that $vx_1\in E(G)$ (by Claim~\ref{claim4.8}(i)). Then
by Lemma~\ref{lem3.2} and Claim~\ref{claim4.12}, there exists a
shadow $u_1$ of $v$ in $V_2$ such that $u_1x_1,u_1x_2\in E(G)$ and
$u_1x_3\notin E(G)$. By Lemma~\ref{lem3.1}(ii), let $u_2$ be a
common neighbor of $v$ and $x_3$. It is clear that
$u_2\notin\{x,x_1,x_2\}$. Then by Claims~\ref{claim4.8}(i)
and~\ref{claim4.12}, we know that $u_2\in V_2\setminus\{u_1\}$. This
shows that $v$ has at least two neighbors in $V_2$, and hence
$f(v)\geq 3$ (by (1)). Suppose $f(v)=3.5$. Then by (1), there must
exist a vertex $v'\in V_1$ such that $vv'\in E(G)$. By the same
argument as above for $v$, we see that $v'$ also has at least two
neighbors in $V_2$ and thus $f(v')=3.5$ (by (1)). But this
contradicts Claim~\ref{claim4.10}. Therefore, we have $f(v)=3$. \qed

\begin{claim}\label{claim4.14}
$G[V_2]$ contains no isolated edges.
\end{claim}

\pf Suppose not, and let $u_1u_2$ be an isolated edge in $G[V_2]$.
Then by Claims~\ref{claim4.8}(i) and~\ref{claim4.12}, we derive that
$u_ix_3\in E(G)$ for some $i\in [2]$; otherwise, $u_1$ and $x_3$ (as
well as $u_2$ and $x_3$) have no common neighbor, contradicting
Lemma~\ref{lem3.1}(ii). By symmetry between $u_1$ and $u_2$ and by
symmetry between $x_1$ and $x_2$, we may assume that
$u_1x_3,u_1x_1\in E(G)$ and $u_1x_2\notin E(G)$.

Suppose $N(u_1)\cap V_1=\emptyset$. Then $N(u_1)=\{x_1,x_3,u_2\}$.
Let $H_1$ be a copy of $W_4$ in $G+u_1x_2$. Since
$x_1x_3,x_2x_3\notin E(G)$, we deduce that $u_1$ is not the center
of $H_1$. This means that the center of $H_1$ is $u_2$, $x_1$ or
$x_2$. But then, it is easy to check that in all cases, there must
exist a vertex $u\in V_2\setminus\{u_1,u_2\}$ such that $uu_2\in
E(G)$, contradicting the assumption that $u_1u_2$ is an isolated
edge in $G[V_2]$.

Thus, we have $N(u_1)\cap V_1\neq\emptyset$. Let $v_1$ be a neighbor
of $u_1$ in $V_1$. Then by Lemma~\ref{lem3.2} and
Claim~\ref{claim4.12}, there exists a shadow $u_3$ of $v_1$ in $V_2$
such that $u_3x_1,u_3x_2\in E(G)$. Notice that $f(v_1)=3$ (by
Claim~\ref{claim4.13}). By (1) and Claim~\ref{claim4.8}(i), we
conclude that $N(v_1)=\{x_j,u_1,u_3\}$ for some $j\in [2]$. Let
$H_2$ be a copy of $W_4$ in $G+v_1x_3$. Then the center of $H_2$ is
$u_1$ (since $u_1$ is the unique common neighbor of $v_1$ and $x_3$
in $G$). Let $v_1x_3y_1y_2v_1$ be the rim of $H_2$. Since
$x_1x_3,x_2x_3\notin E(G)$ and by Claims~\ref{claim4.8}(i)
and~\ref{claim4.12}, we have $y_1\in V_2$. Then, it follows from
$u_1y_1\in E(G)$ and $u_1u_2$ is an isolated edge in $G[V_2]$ that
$y_1=u_2$. Since $v_1y_2,u_2y_2\in E(G)$ and $u_2u_3\notin E(G)$, we
know that $y_2\in\{x_1,x_2\}$. Moreover, because $u_1y_2,u_1x_1\in
E(G)$ and $u_1x_2\notin E(G)$, we see that $y_2=x_1$. This implies
that $N(v_1)=\{x_1,u_1,u_3\}$ and $u_2x_1,u_2x_3\in E(G)$.

Note that neither $u_1$ nor $u_2$ is the shadow of the vertices in
$V_1$ (by Lemma~\ref{lem3.2}). If there exists a vertex $v\in V_1$
such that $vu_1,vu_2\in E(G)$, then by Lemma~\ref{lem3.2} and
Claim~\ref{claim4.12}, there must exist a shadow of $v$ in
$V_2\setminus\{u_1,u_2\}$ and thus $f(v)\geq 4$ (by (1)),
contradicting Claim~\ref{claim4.13}. Hence, $u_1$ and $u_2$ have no
common neighbor in $V_1$. Let $H_3$ be a copy of $W_4$ in
$G+v_1u_2$. Since $u_1u_3,u_2u_3\notin E(G)$, we derive that $v_1$
is not the center of $H_3$. If $u_1$ or $u_2$ is the center of
$H_3$, then it is easy to see that there must exist a vertex $v'\in
V_1$ such that $v'x_1,v'u_1,v'u_2\in E(G)$ (since $u_1u_2$ is an
isolated edge in $G[V_2]$), contradicting the fact that $u_1$ and
$u_2$ have no common neighbor in $V_1$. Therefore, we deduce that
the center of $H_3$ is $x_1$. Let $v_1u_2v_2zv_1$ be the rim of
$H_3$. Since $v_1z\in E(G)$, we have $z\in\{u_1,u_3\}$. If $z=u_1$,
then $v_2x_1,v_2u_1,v_2u_2\in E(G)$, which means that $v_2$ is a
common neighbor of $u_1$ and $u_2$ in $V_1$ (since $u_1u_2$ is an
isolated edge in $G[V_2]$), a contradiction. Thus, $z=u_3$ and
$v_2x_1,v_2u_2,v_2u_3\in E(G)$. Since $u_1u_2$ is an isolated edge
in $G[V_2]$, we can conclude that $v_2\in V_1$. Then by (1) and
Claim~\ref{claim4.13}, we have $N(v_2)=\{x_1,u_2,u_3\}$.

\vspace{15pt}
\begin{figure}[ht]
\begin{center}
\includegraphics[scale=0.8]{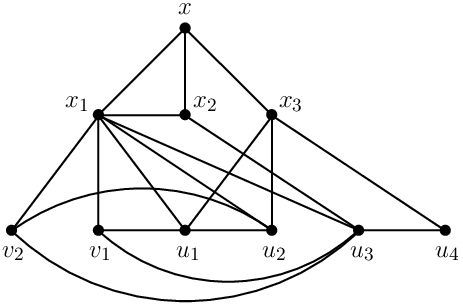}
\caption{The configuration in the proof of Claim~\ref{claim4.14}.}
\label{fig6}
\end{center}
\end{figure}

Since $f(v_1)=f(v_2)=3$ (by Claim~\ref{claim4.13}), it follows from
Claim~\ref{claim4.10} that $V_1=V_x^*=\{v_1,v_2\}$ and hence
$f(u)=2.5$ for any $u\in V_2$. Then $N(u_1)=\{x_1,x_3,u_2,v_1\}$,
$N(u_2)=\{x_1,x_3,u_1,v_2\}$ and $G[V_2]$ is a matching (by (1)).
Let $u_4$ be the unique neighbor of $u_3$ in $V_2$. Then
$N(u_3)=\{x_1,x_2,u_4,v_1,v_2\}$. Moreover, $u_4x_3\in E(G)$;
otherwise, $u_3$ and $x_3$ have no common neighbor, contradicting
Lemma~\ref{lem3.1}(ii). See Figure~\ref{fig6} for an illustration.

Let $H_4$ be a copy of $W_4$ in $G+v_2x_2$. Since
$u_2x_2,u_2u_3\notin E(G)$, we know that $v_2$ is not the center of
$H_4$. This shows that the center of $H_4$ is $u_3$, $x_1$ or $x_2$.
Then, it is straightforward to verify that in all cases, we always
have $V(H_4)=\{u_3,u_4,v_2,x_1,x_2\}$ and $u_4x_1,u_4x_2\in E(G)$.
But this contradicts the fact that $u_4\in V_2$ and $u_4x_3\in
E(G)$. \qed

\begin{claim}\label{claim4.15}
$G[V_2]\in\{P_3,P_4,S_4,2P_3\}$, where $2P_3$ denotes the disjoint
union of $2$ copies of $P_3$.
\end{claim}

\pf Recall that $2.5\leq f(v)\leq 3.5$ for any $v\in V_x$. By (1)
and Claim~\ref{claim4.12}, we see that for any $u\in V_2$,
$f(u)=2.5$ if and only if $u$ has exactly one neighbor in $V_2$,
$f(u)=3$ if and only if $u$ has exactly two neighbors in $V_2$ and
$f(u)=3.5$ if and only if $u$ has exactly three neighbors in $V_2$.
Define $V_2^*:=\{u\in V_2: 3\leq f(u)\leq 3.5\}$. By
Claims~\ref{claim4.8}(ii),~\ref{claim4.12} and~\ref{claim4.14}, we
derive that every component of $G[V_2]$ contains at least three
vertices and thus $|V_2^*|\geq 1$. On the other hand, if follows
from Claim~\ref{claim4.10} that $|V_2^*|\leq |V_x^*|\leq 2$. Then
$1\leq |V_2^*|\leq 2$.

First, suppose $|V_2^*|=1$. Then $G[V_2]$ contains exactly one
component. Let $V_2^*=\{u\}$. It is easy to observe that $G[V_2]$ is
isomorphic to $P_3$ (if $f(u)=3$) or $S_4$ (if $f(u)=3.5$).

Next, suppose $|V_2^*|=2$. Let $V_2^*=\{u_1,u_2\}$. Since
$V_2^*\subseteq V_x^*$ and by Claim~\ref{claim4.10}, we have
$f(u_1)=f(u_2)=3$. If $u_1$ and $u_2$ are contained in the same
component of $G[V_2]$, then $G[V_2]$ is connected and isomorphic to
$P_4$. If $u_1$ and $u_2$ are contained in different components of
$G[V_2]$, then we can deduce that $G[V_2]$ contains exactly two
components, each of which is isomorphic to $P_3$ (i.e. $G[V_2]\cong
2P_3$).   \qed

\bigskip

We now consider two cases according to whether $V_1=\emptyset$ or
not.

\medskip

{\bf Case 1.} $V_1\neq\emptyset$.

\medskip

By Claims~\ref{claim4.11} and~\ref{claim4.15}, we conclude that
$V_x^*\cap V_1\neq\emptyset$ and $V_x^*\cap V_2\neq\emptyset$. Then
by Claim~\ref{claim4.10}, we have $|V_x^*\cap V_1|=|V_x^*\cap
V_2|=1$ and $f(u)=3$ for any $u\in V_x^*$. This implies that
$|V_1|=1$ (since $V_1\subseteq V_x^*$ by Claim~\ref{claim4.13}) and
$G[V_2]\cong P_3$ (by Claim~\ref{claim4.15}). Let $V_1=\{v\}$. Then
by (1), we know that $d(v)=3$ and thus $v$ has two neighbors in
$V_2$. Moreover, $vx_3\notin E(G)$ by Claim~\ref{claim4.8}(i).
Hence, we may assume without loss of generality that
$N(v)=\{x_1,u_1,u_2\}$, where $u_1$ is a shadow of $v_1$ in $V_2$
and $u_1x_1,u_1x_2\in E(G)$ (by Lemma~\ref{lem3.2}). Then by
Lemma~\ref{lem3.1}(ii), we see that $u_2x_3\in E(G)$; otherwise, $v$
and $x_3$ have no common neighbor.

Let $H$ be a copy of $W_4$ in $G+vx_3$. Then the center of $H$ must
be $u_2$ (since $u_2$ is the unique common neighbor of $v$ and $x_3$
in $G$). Let $vx_3u_3yv$ be the rim of $H$. Since $vy,u_3x_3\in
E(G)$ and $x_1x_3,x_2x_3\notin E(G)$, we have $y\in\{x_1,u_1\}$ and
$u_3\in V_2$. If $y=u_1$, then we can derive that
$u_1u_2,u_1u_3,u_2u_3\in E(G)$, contradicting the fact that
$G[V_2]\cong P_3$. Therefore, $y=x_1$. This shows that
$u_2x_1,u_2u_3,u_3x_1,u_3x_3\in E(G)$. Since $G[V_2]\cong P_3$, we
deduce that $V_2=\{u_1,u_2,u_3\}$ and $u_1$ is adjacent to exactly
one vertex in $\{u_2,u_3\}$. But then, one can easily check that
$G+x_2x_3$ contains no copy of $W_4$ (since $x$ is the unique common
neighbor of $x_2$ and $x_3$ in $G$ and $d(x)=3$), a contradiction.

\medskip

{\bf Case 2.} $V_1=\emptyset$.

\medskip

By Claim~\ref{claim4.15}, we notice that $G[V_2]\in\{P_3,P_4,S_4,2P_3\}$.

\medskip

{\bf Subcase 2.1.} $G[V_2]\in\{P_3,2P_3\}$.

\medskip

Let $u_1u_2u_3$ be a copy of $P_3$ in $G[V_2]$. Then $d(u_2)=4$ and
$d(u_1)=d(u_3)=3$. It follows from Claim~\ref{claim4.8}(iii) that
$u_1x_3,u_3x_3\in E(G)$. By symmetry between $x_1$ and $x_2$, we may
also assume that $u_1x_1\in E(G)$. Let $H$ be a copy of $W_4$ in
$G+u_1x_2$. Since $x_1x_3,x_2x_3\notin E(G)$, we conclude that $u_1$
is not the center of $H$. If $x_1$ or $x_2$ is the center of $H$,
then it is easy to see that $V(H)=\{u_1,u_2,u_3,x_1,x_2\}$ and
$u_3x_1,u_3x_2\in E(G)$, contradicting the fact that $u_3\in V_2$
and $u_3x_3\in E(G)$. Thus, we know that the center of $H$ must be
$u_2$ and the rim of $H$ must be $u_1x_2u_3x_3u_1$. This means that
$u_2x_2,u_2x_3,u_3x_2\in E(G)$. Now, it is easy to check that
$G':=G+u_1u_3$ contains no copy of $W_4$ (since $G'[N(u)]$ contains
no copy of $C_4$ for any $u\in\{u_1,u_2,u_3,x_3\}$), giving a
contradiction.

\medskip

{\bf Subcase 2.2.} $G[V_2]\cong S_4$.

\medskip

In this subcase, we apply a similar argument to that in the proof of
Subcase 2.1. Let $V_2=\{u_1,u_2,u_3,u_4\}$ such that
$u_1u_2,u_2u_3,u_2u_4\in E(G)$. Then $d(u_2)=5$ and
$d(u_1)=d(u_3)=d(u_4)=3$. By Claim~\ref{claim4.8}(iii), we have
$u_1x_3,u_3x_3,u_4x_3\in E(G)$. By symmetry between $x_1$ and $x_2$,
we may further assume that $u_1x_1\in E(G)$. Let $H$ be a copy of
$W_4$ in $G+u_1x_2$. Since $x_1x_3,x_2x_3\notin E(G)$, we see that
$u_1$ is not the center of $H$. If $x_1$ or $x_2$ is the center of
$H$, then there must exist some $i\in\{3,4\}$ such that
$V(H)=\{u_1,u_2,u_i,x_1,x_2\}$ and $u_ix_1,u_ix_2\in E(G)$, which
contradicts the fact that $u_i\in V_2$ and $u_ix_3\in E(G)$. Hence,
we can derive that the center of $H$ must be $u_2$ and the rim of
$H$ must be $u_1x_2u_ix_3u_1$ for some $i\in \{3,4\}$. This implies
that $u_2x_2,u_2x_3,u_ix_2\in E(G)$. But then, it is straightforward
to verify that $G':=G+u_1u_i$ contains no copy of $W_4$ (since
$G'[N(u)]$ contains no copy of $C_4$ for any
$u\in\{u_1,u_2,u_i,x_3\}$), a contradiction.

\medskip

{\bf Subcase 2.3.} $G[V_2]\cong P_4$.

\medskip

Let $V_2=\{u_1,u_2,u_3,u_4\}$ such that $u_1u_2,u_2u_3,u_3u_4\in
E(G)$. Then $d(u_2)=d(u_3)=4$ and $d(u_1)=d(u_4)=3$. By
Claim~\ref{claim4.8}(iii), we deduce that $u_1x_3,u_4x_3\in E(G)$.
By symmetry between $x_1$ and $x_2$, we may also assume that
$u_1x_1\in E(G)$.

First, suppose $u_2x_1,u_2x_2\in E(G)$ and $u_2x_3\notin E(G)$. If
$u_3x_3\notin E(G)$, then we conclude that $d(x_3)=3$ and
$e(N[x_3])=3<4=e(N[x])$, which contradicts the choice of $x$.
Therefore, $u_3x_3\in E(G)$. Let $H_1$ be a copy of $W_4$ in
$G+u_1x_2$. Since $x_1x_3,x_2x_3\notin E(G)$, we observe that $u_1$
is not the center of $H_1$. This shows that the center of $H_1$ is
$u_2$, $x_1$ or $x_2$. Then one can easily see that in all cases, we
always have $V(H_1)=\{u_1,u_2,u_3,x_1,x_2\}$ and $u_3x_1,u_3x_2\in
E(G)$, contradicting the fact that $u_3\in V_2$ and $u_3x_3\in
E(G)$.

Next, suppose $u_2x_1,u_2x_3\in E(G)$ and $u_2x_2\notin E(G)$. Let
$H_2$ be a copy of $W_4$ in $G+x_2x_3$. Since $u_2u_4\notin E(G)$,
we know that $u_3$ is not the center of $H_2$. If $x_2$ is the
center of $H_2$, then $u_1,u_2\notin V(H_2)$ (since
$u_1x_2,u_2x_2\notin E(G)$) and it is easy to check that there must
exist some $i\in\{3,4\}$ such that $u_ix_1,u_ix_2,u_ix_3\in E(G)$
(no matter $x\in V(H_2)$ or not), contradicting the fact that
$u_i\in V_2$. This implies that the center of $H_2$ must be $x_3$.
Since $x_1x_3\notin E(G)$, we have $x,x_1\notin V(H_2)$ and thus
$|V(H_2)\cap V_2|=3$. But this is impossible since it is easy to
observe that $G[V_2\cup\{x_2\}]$ contains no copy of $C_4$, a
contradiction.

Finally, suppose $u_2x_2,u_2x_3\in E(G)$ and $u_2x_1\notin E(G)$.
But now, since $x_1x_3,x_2x_3,u_2x_1\notin E(G)$, it is
straightforward to check that $G+u_1x$ contains no copy of $W_4$, a
contradiction.

\medskip

To conclude, we derive a contradiction in all cases and hence no
extremal graph exists in this part.   \qed

\subsection{$\delta(G)=3$ and $e(N[x])=5$}

\medskip

In this part, suppose without loss of generality that
$x_1x_2,x_2x_3\in E(G)$ and $x_1x_3\notin E(G)$. Then
$V_3=\emptyset$; otherwise, $G[\{w,x,x_1,x_2,x_3\}]$ contains a copy
of $W_4$ for any $w\in V_3$, a contradiction. This shows that
$V_2\neq\emptyset$. (If $V_1=\emptyset$, then it follows from $n\geq
6$ that $V_2\neq\emptyset$. If $V_1\neq\emptyset$, then by
Lemma~\ref{lem3.2}, we see that $V_2\neq\emptyset$.)

For any $v\in V_1$, we define
\begin{itemize}
\item $R_v:=N(v)\cap V_1$;
\item $S_v:=\{u\in N(v)\cap V_2: N(u)\cap V_2=\emptyset\}$;
\item $T_v:=\{u\in N(v)\cap V_2: N(u)\cap V_2\neq\emptyset\}$.
\end{itemize}
Let $r_v:=|R_v|$, $s_v:=|S_v|$ and $t_v:=|T_v|$, and we say that $v$
is of \emph{Type $(r_v, s_v, t_v)$}. It is clear that for any $v\in
V_1$, we always have $r_v+s_v+t_v=d(v)-1\geq 2$ (since
$V_3=\emptyset$ and $\delta(G)=3$).

\begin{claim}\label{claim4.16}
For any $v\in V_1$ and any shadow $u$ of $v$ in $V_2$, we have $u\in
T_v$.
\end{claim}

\pf Let $v$ be an arbitrary vertex in $V_1$ and let $u$ be any
shadow of $v$ in $V_2$. Since $x_1x_3\notin E(G)$ and by
Lemma~\ref{lem3.2}, we derive that either $ux_1,ux_2\in E(G)$ or
$ux_2,ux_3\in E(G)$. If $N(u)\cap V_2=\emptyset$, then by
Lemma~\ref{lem3.3}, we deduce that $x_1x_3\in E(G)$ (in both cases),
contradicting the assumption that $x_1x_3\notin E(G)$. Thus, we have
$N(u)\cap V_2\neq\emptyset$, i.e. $u\in T_v$.     \qed

\bigskip

In this and the next subsection, we shall use the discharging
method. For any $v\in V_1\cup V_2$, let $f(v)$ be the initial charge
of $v$. Then we redistribute the charges of the vertices in $V_1\cup
V_2$ according to the following discharging rule:
\begin{itemize}
\item [(\bf R)] For any $v\in V_1$, if $0.5r_v+0.5s_v+t_v\geq 1.5$, then
$v$ sends $0.5$ to each of its neighbors in $S_v$.
\end{itemize}
Let $f^*(v)$ be the new charge of $v$ for any $v\in V_1\cup V_2$
after applying the above discharging rule. Since $V_3=\emptyset$, it
is obvious that
\begin{align*}
\sum\limits_{v\in V_x}f^*(v)=\sum\limits_{v\in V_x}f(v).
\end{align*}
Then by Lemma~\ref{lem3.4}, we conclude that
\begin{align*}
e(G)=e(N[x])+\sum\limits_{v\in V_x}f^*(v). \tag{$3$}
\end{align*}

\begin{claim}\label{claim4.17}
$f^*(v)\geq 2.5$ for any $v\in V_1\cup V_2$.
\end{claim}

\pf First, suppose $v\in V_1$. Then by Lemma~\ref{lem3.2}, we know
that $v$ has at least one shadow in $V_2$. This, together with
Claim~\ref{claim4.16}, implies that $t_v\geq 1$. Since
$r_v+s_v+t_v\geq 2$, we have
\begin{align*}
0.5r_v+0.5s_v+t_v=0.5(r_v+s_v+t_v)+0.5t_v\geq 0.5\cdot 2+0.5\cdot
1=1.5
\end{align*}
with equality if and only if $r_v=0$ and $s_v=t_v=1$, or $s_v=0$ and
$r_v=t_v=1$ (i.e. $v$ is of Type $(0,1,1)$ or
Type $(1,0,1)$). Then by the discharging rule (R), we see that $v$ sends
$0.5$ to each of its neighbors in $S_v$. Hence, it follows from (1)
that
\begin{align*}
f^*(v)=f(v)-0.5s_v=1+0.5r_v+(s_v+t_v)-0.5s_v=1+(0.5r_v+0.5s_v+t_v)\geq
2.5,
\end{align*}
and the equality holds if and only if $v$ is of Type $(0,1,1)$ or
Type $(1,0,1)$.

Next, suppose $v\in V_2$. If $N(v)\cap V_2\neq\emptyset$, then
$v\notin S_{v'}$ for any $v'\in V_1$ (by the definition of $S_{v'}$)
and thus $f^*(v)=f(v)\geq 2.5$ with equality if and only if
$|N(v)\cap V_2|=1$ (by (1)). So we may assume that $N(v)\cap
V_2=\emptyset$. Since $\delta(G)=3$, there exists a vertex $v'\in
V_1$ such that $vv'\in E(G)$. By the same argument as above for $v$,
we can show that $0.5r_{v'}+0.5s_{v'}+t_{v'}\geq 1.5$. Then by the
discharging rule (R), we derive that $v$ receives $0.5$ from $v'$.
Therefore,
\begin{align*}
f^*(v)\geq f(v)+0.5=2+0.5=2.5
\end{align*}
(by (1)), and the equality holds if and only if $v'$ is the unique
neighbor of $v$ in $V_1$ (i.e. $|N(v)\cap V_1|=1$).   \qed

\bigskip

Now, by (3) and Claim~\ref{claim4.17}, we deduce that
\begin{align*}
e(G)\geq
5+2.5(n-4)=\frac{5n-10}{2}\geq\lfloor\frac{5n-10}{2}\rfloor.
\end{align*}

In the following, we characterize the extremal graphs. Suppose
$e(G)=\lfloor\frac{5n-10}{2}\rfloor$. If there exists some vertex
$v\in V_1\cup V_2$ such that $f^*(v)\geq 3$, then by (3) and
Claim~\ref{claim4.17}, we have
\begin{align*}
e(G)\geq 5+3+2.5(n-5)=\frac{5n-9}{2}>\lfloor\frac{5n-10}{2}\rfloor,
\end{align*}
a contradiction. Thus, we conclude that $f^*(v)=2.5$ for
any $v\in V_1\cup V_2$. Then, it follows from the proof of
Claim~\ref{claim4.17} that the following statements hold:
\begin{itemize}
\item [(Q1)] if $v\in V_1$, then $v$ is of Type $(0,1,1)$ or
Type $(1,0,1)$;
\item [(Q2)] if $v\in V_2$ and $N(v)\cap V_2\neq\emptyset$, then
$|N(v)\cap V_2|=1$;
\item [(Q3)] if $v\in V_2$ and $N(v)\cap V_2=\emptyset$, then
$|N(v)\cap V_1|=1$.
\end{itemize}

\begin{claim}\label{claim4.18}
Every vertex in $V_1$ is of Type $(1,0,1)$.
\end{claim}

\pf Suppose this is false. Then by (Q1), there must exist a vertex
$v\in V_1$ such that $v$ is of Type $(0,1,1)$. Let $u_1$ be the
unique neighbor of $v$ in $S_v$ and $u_2$ the unique neighbor of $v$
in $T_v$. Then by Claim~\ref{claim4.16}, we know that $u_2$ is the
unique shadow of $v$ in $V_2$. By (Q2), let $u_3$ be the unique
neighbor of $u_2$ in $V_2$. If $u_1x_1,u_1x_2\in E(G)$ or
$u_1x_2,u_1x_3\in E(G)$, then by Lemma~\ref{lem3.3}, we have
$x_1x_3\in E(G)$, which contradicts the assumption that
$x_1x_3\notin E(G)$. Hence, we see that $u_1x_1,u_1x_3\in E(G)$, and
thus $N(u_1)=\{x_1,x_3,v\}$ (by (Q3)).

Suppose $vx_2\notin E(G)$. Then by symmetry between $x_1$ and $x_3$,
we may assume that $vx_1\in E(G)$. Since $x_1x_3\notin E(G)$ and by
Lemma~\ref{lem3.2}, we derive that $u_2x_1,u_2x_2\in E(G)$. But
then, one can easily check that $G+vx_3$ contains no copy of $W_4$
(since $u_1$ is the unique common neighbor of $v$ and $x_3$ in $G$
and $d(u_1)=3$), a contradiction.

Therefore, we have $vx_2\in E(G)$, and hence $N(v)=\{x_2,u_1,u_2\}$.
By Lemma~\ref{lem3.2}, we may assume by symmetry that
$u_2x_1,u_2x_2\in E(G)$. Let $H_1$ be a copy of $W_4$ in $G+vx_3$.
Since $u_1u_2\notin E(G)$, we deduce that $v$ is not the center of
$H_1$. Moreover, because $x_1x_3,u_2x_3\notin E(G)$, we observe that
$x_3$ is also not the center of $H_1$. Then, it is easy to see that
the center of $H_1$ must be $x_2$ and the rim of $H_1$ must be
$vx_3u_3u_2v$. This shows that $u_3x_2,u_3x_3\in E(G)$.

Let $H_2$ be a copy of $W_4$ in $G+vx_1$. Since $u_1u_2,u_1x_2\notin
E(G)$, we conclude that $v$ is not the center of $H_2$. On the other
hand, because $u_3v,u_3x_1\notin E(G)$, we know that $u_2$ is also
not the center of $H_2$. This implies that the center of $H_2$ is
$x_1$ or $x_2$. Since $x_1x_3\notin E(G)$, it is straightforward to
verify that in both cases, we always have
$V(H_2)=\{v,u_2,u_3,x_1,x_2\}$ and $u_3x_1\in E(G)$. But this
contradicts the fact that $u_3\in V_2$ and $u_3x_2,u_3x_3\in E(G)$.
\qed

\begin{claim}\label{claim4.19}
If $u_1u_2$ is an edge in $G[V_2]$ such that $d(u_1)=3$ and
$u_1x_2\notin E(G)$, then $u_2x_2\notin E(G)$.
\end{claim}

\pf Suppose to the contrary that $u_2x_2\in E(G)$. Then $u_2$ is
adjacent to exactly one vertex in $\{x_1,x_3\}$ (because $u_2\in
V_2$). Since $d(u_1)=3$ and $u_1x_2\notin E(G)$, we have
$N(u_1)=\{x_1,x_3,u_2\}$. But then, we see that $d(u_1)=3$ and
$e(N[u_1])=4<5=e(N[x])$, which contradicts the choice of $x$. \qed

\bigskip

We now consider two cases according to whether $V_1=\emptyset$ or
not.

\medskip

{\bf Case 1.} $V_1\neq\emptyset$.

\medskip

Suppose there exists a vertex $v_1\in V_1$ such that $v_1x_2\notin
E(G)$. Then by symmetry, we may assume that $v_1x_1\in E(G)$. Since
$v_1$ is of Type $(1,0,1)$ (by Claim~\ref{claim4.18}), we may
further assume that $N(v_1)=\{x_1,v_2,u\}$, where $v_2$ is the
unique neighbor of $v_1$ in $V_1$ and $u$ is the unique neighbor of
$v_1$ in $V_2$. Then by Lemma~\ref{lem3.2}, we derive that
$ux_1,ux_2\in E(G)$. This means that $v_2x_3\in E(G)$; otherwise,
$v_1$ and $x_3$ have no common neighbor, contradicting
Lemma~\ref{lem3.1}(ii). Moreover, $v_2u\in E(G)$; otherwise, we
deduce that $d(v_1)=3$ and $e(N[v_1])=4<5=e(N[x])$, which
contradicts the choice of $x$. Since $v_2$ is also of Type $(1,0,1)$
(by Claim~\ref{claim4.18}), we have $N(v_2)=\{x_3,v_1,u\}$ and hence
$u$ is the unique shadow of $v_2$ in $V_2$. But this implies that
$ux_2,ux_3\in E(G)$ (by Lemma~\ref{lem3.2}), contradicting the fact
that $u\in V_2$ and $ux_1,ux_2\in E(G)$.

Thus, we conclude that $vx_2\in E(G)$ for any $v\in V_1$. Since $G$
contains no universal vertex, we know that there must exist a vertex
$u_1\in V_2$ such that $u_1x_2\notin E(G)$; otherwise, $x_2$ would
be a universal vertex of $G$, a contradiction. If there exists some
vertex $v\in V_1$ such that $vu_1\in E(G)$, then it follows from
Claim~\ref{claim4.18} and Lemma~\ref{lem3.2} that $u_1$ is the
unique shadow of $v$ in $V_2$ and $u_1x_2\in E(G)$ (since $vx_2\in
E(G)$), a contradiction. Hence, we see that $N(u_1)\cap
V_1=\emptyset$, and thus $N(u_1)\cap V_2\neq\emptyset$ (since
$\delta(G)=3$). Then by (Q2), we may assume that
$N(u_1)=\{x_1,x_3,u_2\}$, where $u_2$ is the unique neighbor of
$u_1$ in $V_2$. Since $u_1x_2\notin E(G)$ and by
Claim~\ref{claim4.19}, we have $u_2x_2\notin E(G)$. Then by the same
argument as above for $u_1$, we can also derive that $N(u_2)\cap
V_1=\emptyset$. But then, $u_1$ and $v$ have no common neighbor for
any $v\in V_1$ (since $vx_2\in E(G)$), contradicting
Lemma~\ref{lem3.1}(ii).

\medskip

{\bf Case 2.} $V_1=\emptyset$.

\medskip

By (Q2) and (Q3), we deduce that every vertex in $V_2$ has exactly
one neighbor in $V_2$. This shows that $d(u)=3$ for any $u\in V_2$
and $G[V_2]$ is a matching. Since $G$ contains no universal vertex,
we conclude that there must exist a vertex $u_1\in V_2$ such that
$u_1x_2\notin E(G)$; otherwise, $x_2$ would be a universal vertex of
$G$, a contradiction. Let $u_2$ be the unique neighbor of $u_1$ in
$V_2$. Then by Claim~\ref{claim4.19}, we have $u_2x_2\notin E(G)$.
Let $H$ be a copy of $W_4$ in $G+x_1x_3$. It is clear that $x\notin
V(H)$ (because $d(x)=3$ and $V_3=\emptyset$). Since $x_1$, $x_2$ and
$x_3$ are the only possible vertices of $G+x_1x_3$ with degree at
least $4$, we know that the center of $H$ must be one vertex in
$\{x_1,x_2,x_3\}$.

First, suppose by symmetry between $x_1$ and $x_3$ that $x_1$ is the
center of $H$. Let $x_3y_1y_2y_3x_3$ be the rim of $H$. Then
$x_2\in\{y_1,y_2,y_3\}$; otherwise, we see that $y_1,y_2,y_3\in
V_2$, which contradicts the fact that $G[V_2]$ is a matching (since
$y_1y_2,y_2y_3\in E(G)$). If $x_2=y_2$, then we can derive that
$y_1x_1,y_1x_2,y_1x_3\in E(G)$ (i.e. $y_1\in V_3$), contradicting
the fact that $V_3=\emptyset$. Therefore, we may assume by symmetry
that $x_2=y_1$. This means that $y_2,y_3\in V_2$ and
$y_2x_1,y_2x_2,y_3x_1,y_3x_3\in E(G)$. But this is impossible since
$d(y_3)=3$ and $y_3x_2\notin E(G)$ would imply that $y_2x_2\notin
E(G)$ (by Claim~\ref{claim4.19}).

Next, suppose $x_2$ is the center of $H$. Let $x_1x_3u_3u_4x_1$ be
the rim of $H$. Then we have $u_3,u_4\in V_2$ and
$u_3u_4,u_3x_2,u_3x_3,u_4x_1,u_4x_2\in E(G)$. Since
$u_1x_2,u_2x_2\notin E(G)$, we deduce that $u_1$, $u_2$, $u_3$ and
$u_4$ are pairwise distinct. But now, it is easy to check that
$G+u_2u_4$ contains no copy of $W_4$ (since $x_1$ is the unique
common neighbor of $u_2$ and $u_4$ in $G$ and $x_1x_3,u_1x_2\notin
E(G)$), giving a contradiction.

\medskip

To sum up, we obtain a contradiction in all cases and thus there is
no extremal graph in this part.   \qed

\subsection{$\delta(G)=3$ and $e(N[x])=6$}

\medskip

In this part, we have $x_1x_2,x_1x_3,x_2x_3\in E(G)$. Moreover, it
follows from the choice of $x$ that $e(N[v])=6$ (i.e. $G[N[v]]\cong
K_4$) for any $v\in V(G)$ with $d(v)=3$. By the same argument as
that in Subsection 4.4, we can conclude that $V_3=\emptyset$ and
$V_2\ne\emptyset$.

For any $v\in V_1$, let $R_v$, $S_v$, $T_v$, $r_v$, $s_v$ and $t_v$
be defined the same as that in Subsection 4.4. Then, it is easy to
observe that $r_v+s_v+t_v=d(v)-1\geq 2$ (since $V_3=\emptyset$ and
$\delta(G)=3$).

\begin{claim}\label{claim4.20}
No vertex in $V_1$ is of Type $(0,2,0)$.
\end{claim}

\pf Suppose not, and let $v$ be a vertex in $V_1$ such that $v$ is
of Type $(0,2,0)$. Let $u_1$ and $u_2$ be the two neighbors of $v$
in $S_v$. Then $u_1u_2\notin E(G)$ (by the definition of $S_v$). But
then, we notice that $d(v)=3$ and $e(N[v])\leq 5<6=e(N[x])$, which
contradicts the choice of $x$.    \qed

\bigskip

Let the discharging rule (R) be defined the same as that in
Subsection 4.4. For any $v\in V_1\cup V_2$, we still let $f(v)$ be
the initial charge of $v$ and $f^*(v)$ the new charge of $v$ after
applying the discharging rule (R). Define $V_2^*:=\{u\in V_2:
N(u)\cap V_2=\emptyset \mbox{ and } v \mbox{ is of Type } (1,1,0)
\mbox{ for any } v\in N(u)\cap V_1\}$.

\begin{claim}\label{claim4.21}
$f^*(v)\geq 2.5$ for any $v\in V_1\cup(V_2\setminus V_2^*)$.
\end{claim}

\pf First, suppose $v\in V_1$. If $v$ does not satisfy the condition
$0.5r_v+0.5s_v+t_v\geq 1.5$, then by (1) and Lemma~\ref{lem3.2}, we
know that $f^*(v)=f(v)\geq 2.5$. So we may assume that $v$ satisfies
the condition $0.5r_v+0.5s_v+t_v\geq 1.5$. Then by the discharging
rule (R), we see that $v$ sends $0.5$ to each of its neighbors in
$S_v$. Thus, it follows from (1) that
\begin{align*}
f^*(v)=f(v)-0.5s_v=1+0.5r_v+(s_v+t_v)-0.5s_v=1+(0.5r_v+0.5s_v+t_v)\geq
2.5.
\end{align*}

Next, suppose $v\in V_2\setminus V_2^*$. If $N(v)\cap
V_2\neq\emptyset$, then $v\notin S_{v'}$ for any $v'\in V_1$ (by the
definition of $S_{v'}$) and thus $f^*(v)=f(v)\geq 2.5$ (by (1)).
Hence, we may assume that $N(v)\cap V_2=\emptyset$. Since
$\delta(G)=3$, we have $N(v)\cap V_1\neq\emptyset$. Then by the
definition of $V_2^*$, there must exist a vertex $v'\in N(v)\cap
V_1$ such that $v'$ is not of Type $(1,1,0)$ (since $v\in
V_2\setminus V_2^*$). By Claim~\ref{claim4.20}, $v'$ is also not of
Type $(0,2,0)$. Note that $r_{v'}+s_{v'}+t_{v'}\geq 2$ and
$s_{v'}\geq 1$ (since $v\in S_{v'}$). Therefore, we can derive that
either $v'$ is of Type $(0,1,1)$ or $r_{v'}+s_{v'}+t_{v'}\geq 3$. In
both cases, we have
\begin{align*}
0.5r_{v'}+0.5s_{v'}+t_{v'}=0.5(r_{v'}+s_{v'}+t_{v'})+0.5t_{v'}\geq
1.5.
\end{align*}
Then by the discharging rule (R), we deduce that $v$ receives $0.5$
from $v'$. By (1), we conclude that
\begin{align*}
f^*(v)\geq f(v)+0.5=2+0.5=2.5.
\end{align*}
This completes the proof of the claim.     \qed

\begin{claim}\label{claim4.22}
$|V_2^*|\leq 1$.
\end{claim}

\begin{figure}[ht]
\begin{center}
\includegraphics[scale=0.8]{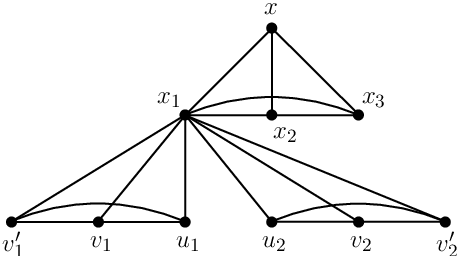}
\caption{The configuration in the proof of Claim~\ref{claim4.22}.}
\label{fig7}
\end{center}
\end{figure}

\pf Suppose not, and let $u_1$ and $u_2$ be two vertices in $V_2^*$.
Then $N(u_1)\cap V_2=N(u_2)\cap V_2=\emptyset$. Since $\delta(G)=3$,
there exists a vertex $v_i\in V_1$ such that $v_iu_i\in E(G)$ for
each $i\in [2]$. By the definition of $V_2^*$, we know that both
$v_1$ and $v_2$ are of Type $(1,1,0)$ (and hence $v_1\neq v_2$). For
each $i\in [2]$, let $v_i'$ be the unique neighbor of $v_i$ in
$V_1$. Without loss of generality, we may assume that
$N(v_1)=\{x_1,v_1',u_1\}$ and $N(v_2)=\{x_i,v_2',u_2\}$ for some
$i\in [3]$. Then by the choice of $x$, we see that
$e(N[v_1])=e(N[v_2])=6$. This implies that
$v_1'u_1,v_1'x_1,u_1x_1,v_2'u_2,v_2'x_i,u_2x_i\in E(G)$. Since
$v_1'\in N(u_1)\cap V_1$ and $v_2'\in N(u_2)\cap V_1$, we can derive
that both $v_1'$ and $v_2'$ are also of Type $(1,1,0)$ (by the
definition of $V_2^*$). Thus, we have $N(v_1')=\{x_1,v_1,u_1\}$ and
$N(v_2')=\{x_i,v_2,u_2\}$. Moreover, $v_1$, $v_1'$, $v_2$ and $v_2'$
are pairwise distinct. This shows that $i=1$; otherwise, $v_1$ and
$v_2$ have no common neighbor, contradicting Lemma~\ref{lem3.1}(ii).
See Figure~\ref{fig7} for an illustration.

\vspace{15pt}

Let $H$ be a copy of $W_4$ in $G+v_1v_2$. Then the center of $H$
must be $x_1$ (since $x_1$ is the unique neighbor of $v_1$ and $v_2$
in $G$). Let $v_1v_2zyv_1$ be the rim of $H$. Because $v_1y,v_2z\in
E(G)$, we deduce that $y\in\{v_1',u_1\}$ and $z\in\{v_2',u_2\}$. But
this is impossible since one can easily see that there is no edge
with one endvertex in $\{v_1',u_1\}$ and the other endvertex in
$\{v_2',u_2\}$, contradicting the assumption that $yz\in E(G)$. \qed

\bigskip

If $V_2^*=\emptyset$, then by (3) and Claim~\ref{claim4.21}, we
conclude that
\begin{align*}
e(G)\geq 6+2.5(n-4)=\frac{5n-8}{2}>\lfloor\frac{5n-10}{2}\rfloor.
\end{align*}
Hence by Claim~\ref{claim4.22}, we may assume that
$|V_2^*|=1$. Let $V_2^*=\{u\}$. By the definition of $V_2^*$, it is
easy to verify that $0.5r_v+0.5s_v+t_v=1<1.5$ for any $v\in N(u)\cap
V_1$ (since every such vertex is of Type $(1,1,0)$). Then by (1), we
know that $f^*(u)=f(u)=2$. Now, it follows from (3) and
Claim~\ref{claim4.21} that
\begin{align*}
e(G)\geq 6+2+2.5(n-5)=\frac{5n-9}{2}>\lfloor\frac{5n-10}{2}\rfloor.
\end{align*}

In both cases, we show that $e(G)>\lfloor\frac{5n-10}{2}\rfloor$ and
no extremal graph exists in this part.   \qed

\subsection{$\delta(G)=4$}

\medskip

In this part, we divide the proof into two subsections according to
the value of $n$.

\subsubsection{$6\leq n\leq 11$}

\medskip

In this subsection, we suppose that $6\leq n\leq 11$.

If there exists some vertex $v\in V(G)$ such that $d(v)\geq 6$ or
two vertices $v,v'\in V(G)$ such that $d(v)=d(v')=5$, then we have
\begin{align*}e(G)\geq
\frac{6+4(n-1)}{2}=2n+1=\frac{5n-9}{2}+\frac{11-n}{2}\geq
\frac{5n-9}{2}>\lfloor\frac{5n-10}{2}\rfloor
\end{align*}
or
\begin{align*}
e(G)\geq\frac{5\cdot
2+4(n-2)}{2}=2n+1=\frac{5n-9}{2}+\frac{11-n}{2}\geq
\frac{5n-9}{2}>\lfloor\frac{5n-10}{2}\rfloor.
\end{align*}
So we may
assume that every vertex in $G$ has degree $4$ or $5$ and the number
of vertices of degree $5$ in $G$ is at most one. Since every graph
contains an even number of vertices of odd degree, we see that $G$
is $4$-regular.

If $6\leq n\leq 9$, then we derive that
\begin{align*}
e(G)=2n=\frac{5n-9}{2}+\frac{9-n}{2}\geq
\frac{5n-9}{2}>\lfloor\frac{5n-10}{2}\rfloor.
\end{align*}
Therefore, we always assume that $G$ is $4$-regular and $10\leq
n\leq 11$ in the rest of this subsection.

\begin{claim}\label{claim4.23}
$|N(u_1)\cap N(u_2)|\geq 2$ for any pair of non-adjacent vertices
$u_1$ and $u_2$ in $G$.
\end{claim}

\pf Suppose not, and let $u_1$ and $u_2$ be two non-adjacent
vertices in $G$ such that $|N(u_1)\cap N(u_2)|\leq 1$. Then by
Lemma~\ref{lem3.1}(ii), we have $|N(u_1)\cap N(u_2)|=1$. Let
$N(u_1)=\{v,y_1,y_2,y_3\}$ and $N(u_2)=\{v,z_1,z_2,z_3\}$, where $v$
is the unique common neighbor of $u_1$ and $u_2$ in $G$. Let $H_1$
be a copy of $W_4$ in $G+u_1u_2$. It is straightforward to check
that the center of $H_1$ must be $v$ (since $v$ is the unique common
neighbor of $u_1$ and $u_2$ in $G$) and the rim of $H_1$ must be
$u_1u_2z_iy_ju_1$ for some $i,j\in [3]$. Without loss of generality,
we may assume that $i=j=1$. This implies that $vy_1,vz_1,y_1z_1\in
E(G)$.

Since $G$ is $4$-regular, there exists some vertex in
$\{y_2,y_3,z_2,z_3\}$, say $y_2$, such that $y_2y_1,y_2z_1\notin
E(G)$. Let $H_2$ be a copy of $W_4$ in $G+vy_2$. Then the center of
$H_2$ must be $u_1$ (since $u_1$ is the unique common neighbor of
$v$ and $y_2$ in $G$) and the rim of $H_2$ must be $vy_2y_3y_1v$.
This shows that $y_3y_1,y_3y_2\in E(G)$. Since $G$ is $4$-regular,
we notice that $z_1$ has at most one neighbor in $\{z_2,z_3\}$. By
symmetry, we may assume that $z_1z_2\notin E(G)$. Then $y_3z_2\in
E(G)$; otherwise, $y_1$ and $z_2$ have no common neighbor,
contradicting Lemma~\ref{lem3.1}(ii).

\begin{figure}[ht]
\begin{center}
\includegraphics[scale=0.8]{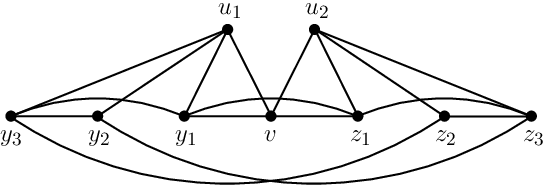}
\caption{The configuration in the proof of Claim~\ref{claim4.23}.}
\label{fig8}
\end{center}
\end{figure}

Let $H_3$ be a copy of $W_4$ in $G+vz_2$. Then, it is easy to see
that the center of $H_3$ must be $u_2$ (since $u_2$ is the unique
common neighbor of $v$ and $z_2$ in $G$) and the rim of $H_3$ must
be $vz_2z_3z_1v$. This means that $z_3z_1,z_3z_2\in E(G)$. By
Lemma~\ref{lem3.1}(ii), we deduce that $z_3y_2\in E(G)$; otherwise,
$z_1$ and $y_2$ have no common neighbor. See Figure~\ref{fig8} for
an illustration.

\vspace{15pt}

Let $U:=\{u_1,u_2,v,y_1,y_2,y_3,z_1,z_2,z_3\}$ and $W:=V(G)\setminus
U$. Then $1\leq |W|\leq 2$ (because $10\leq n\leq 11$). Since $G$ is
$4$-regular, we conclude that $e(U,W)\leq 2$. But this implies that
$d(w)\leq 3$ for any $w\in W$, a contradiction.  \qed

\begin{claim}\label{claim4.24}
$|N(u_1)\cap N(u_2)|\neq 3$ for any pair of non-adjacent vertices
$u_1$ and $u_2$ in $G$.
\end{claim}

\pf Suppose not, and let $u_1$ and $u_2$ be two non-adjacent
vertices in $G$ such that $|N(u_1)\cap N(u_2)|=3$. Let
$N(u_1)=\{v_1,v_2,v_3,y\}$ and $N(u_2)=\{v_1,v_2,v_3,z\}$, where
$v_1$, $v_2$ and $v_3$ are the three common neighbors of $u_1$ and
$u_2$ in $G$. Define $W:=V(G)\setminus\{u_1,u_2,v_1,v_2,v_3,y,z\}$.
Then $3\leq |W|\leq 4$ (because $10\leq n\leq 11$). Since $|N(w)\cap
N(u_1)|\geq 2$ for any $w\in W$ (by Claim~\ref{claim4.23}), we have
$e(W,N(u_1))\geq 6$.

If $y$ has at least two neighbors in $\{v_1,v_2,v_3\}$, then we can
know that $e(W,N(u_1))\leq 5$ (since $G$ is $4$-regular), a
contradiction. Thus, we see that $y$ has at most one neighbor in
$\{v_1,v_2,v_3\}$. Then, it follows from $|N(y)\cap N(u_2)|\geq 2$
(by Claim~\ref{claim4.23}) that $y$ has exactly one neighbor in
$\{v_1,v_2,v_3\}$ and $yz\in E(G)$. Moreover, because $|N(z)\cap
N(u_1)|\geq 2$ (by Claim~\ref{claim4.23}), we derive that $z$ also
has at least one neighbor in $\{v_1,v_2,v_3\}$. But now, since $G$
is $4$-regular, it is easy to verify that $e(W,N(u_1))\leq 5$,
giving a contradiction. \qed

\begin{claim}\label{claim4.25}
$|N(u_1)\cap N(u_2)|=2$ for any pair of non-adjacent vertices $u_1$
and $u_2$ in $G$.
\end{claim}

\pf Suppose not, and let $u_1$ and $u_2$ be two non-adjacent
vertices in $G$ such that $|N(u_1)\cap N(u_2)|\neq 2$. Then by
Claims~\ref{claim4.23} and~\ref{claim4.24}, we conclude that
$|N(u_1)\cap N(u_2)|=4$. Let $N(u_1)=N(u_2)=\{v_1,v_2,v_3,v_4\}$.
Define $W:=V(G)\setminus\{u_1,u_2,v_1,v_2,v_3,v_4\}$. Then $4\leq
|W|\leq 5$ (since $10\leq n\leq 11$). Because $|N(w)\cap N(u_1)|\geq
2$ for any $w\in W$ (by Claim~\ref{claim4.23}), we have
$e(W,N(u_1))\geq 8$. On the other hand, since $G$ is $4$-regular and
$v_iu_1,v_iu_2\in E(G)$ for each $i\in [4]$, we know that
$e(W,N(u_1))\leq 8$. This shows that $e(W,N(u_1))=8$, and hence
$v_iv_j\notin E(G)$ for any $i,j\in [4]$.

Let $H$ be a copy of $W_4$ in $G+u_1u_2$. It is clear that the
center of $H$ is $u_1$, $u_2$ or $v_i$ for some $i\in [4]$. But
then, one can easily check that in all cases, there must exist some
$p,q,r\in [4]$ such that $v_pv_q,v_qv_r\in E(G)$, contradicting the
fact that $v_iv_j\notin E(G)$ for any $i,j\in [4]$.    \qed

\bigskip

By Claim~\ref{claim4.25}, we observe that every pair of non-adjacent
vertices in $G$ have exactly two common neighbors. We shall use this
fact frequently in the following argument.

Let $u_1$ and $u_2$ be two non-adjacent vertices such that
$N(u_1)=\{v_1,v_2,y_1,y_2\}$ and $N(u_2)=\{v_1,v_2,z_1,z_2\}$, where
$v_1$ and $v_2$ are the two common neighbors of $u_1$ and $u_2$ in
$G$. Let $W:=V(G)\setminus\{u_1,u_2,v_1,v_2,y_1,y_2,z_1,z_2\}$.
Since $10\leq n\leq 11$, we have $2\leq |W|\leq 3$. For the sake of
convenience, we may assume that $W=\{w_1,\ldots,w_k\}$ for some
$k\in\{2,3\}$.

Note that $|N(w_i)\cap N(u_1)|=|N(w_i)\cap N(u_2)|=2$ for each $i\in
[k]$. If $e(W,\{v_1,v_2\})=0$, then we see that
$w_iy_1,w_iy_2,w_iz_1,w_iz_2\in E(G)$ for each $i\in [k]$ and thus
$|N(y_1)\cap N(u_2)|\leq 1$ (since $G$ is $4$-regular), a
contradiction. Hence, we may assume without loss of generality that
$w_1v_1\in E(G)$. Then $w_1v_2\notin E(G)$; otherwise, we have
$|N(v_1)\cap N(v_2)|\geq 3$ and it follows from
Claim~\ref{claim4.25} that $v_1v_2\in E(G)$, which means that
$|N(w_2)\cap N(v_1)|\leq 1$ (since $G$ is $4$-regular), a
contradiction. Because $|N(w_1)\cap N(u_1)|=|N(w_1)\cap N(u_2)|=2$,
we may assume by symmetry that $w_1y_1,w_1z_1\in E(G)$ and
$w_1y_2,w_1z_2\notin E(G)$. Moreover, since $G$ is $4$-regular, we
may further assume that $w_2$ is the unique neighbor of $w_1$ in
$W$. Then $N(w_1)=\{v_1,y_1,z_1,w_2\}$. This implies that
$v_1v_2\notin E(G)$; otherwise, we can derive that $v_2y_2\in E(G)$
(since $|N(y_2)\cap N(v_1)|=2$) and thus $|N(w_2)\cap N(v_1)|=1$
(since $G$ is $4$-regular), giving a contradiction.

First, we consider the vertex $w_2$. Suppose $w_2v_1\in E(G)$. Since
$|N(y_2)\cap N(v_1)|=|N(z_2)\cap N(v_1)|=2$, we can deduce that
$w_2y_2,w_2z_2\in E(G)$. This shows that $v_2y_2,v_2z_2\in E(G)$
(because $|N(w_2)\cap N(v_2)|=2$). But then, we notice that
$N(w_1)\cap N(v_2)=\emptyset$ (since $N(w_1)=\{v_1,y_1,z_1,w_2\}$
and $N(v_2)=\{u_1,u_2,y_2,z_2\}$), a contradiction. Therefore,
$w_2v_1\notin E(G)$. Then, we have $w_2v_2\in E(G)$; otherwise, we
conclude that $w_2y_1,w_2y_2\in E(G)$ (since $|N(w_2)\cap
N(u_1)|=2$), which means that $|N(w_2)\cap N(u_2)|\leq 1$ (since $G$
is $4$-regular), a contradiction.

Next, we consider the vertex $y_2$. Suppose $y_2v_1\in E(G)$. Then,
we know that $y_2w_2,y_2z_2\in E(G)$ (because $|N(w_2)\cap
N(v_1)|=|N(z_2)\cap N(v_1)|=2$). Since $|N(y_2)\cap N(z_1)|=2$, we
have $z_1w_2,z_1z_2\in E(G)$. But this implies that $N(z_1)\cap
N(u_1)=\emptyset$ (because $N(z_1)=\{u_2,z_2,w_1,w_2\}$ and
$N(u_1)=\{v_1,v_2,y_1,y_2\}$), a contradiction. Thus, we see that
$y_2v_1\notin E(G)$. Moreover, $y_2v_2\notin E(G)$; otherwise, we
derive that $|N(w_1)\cap N(v_2)|=1$ (since $G$ is $4$-regular),
giving a contradiction. This shows that $y_2z_1,y_2z_2\in E(G)$
(because $|N(y_2)\cap N(u_2)|=2$).

\vspace{15pt}
\begin{figure}[ht]
\begin{center}
\includegraphics[scale=0.8]{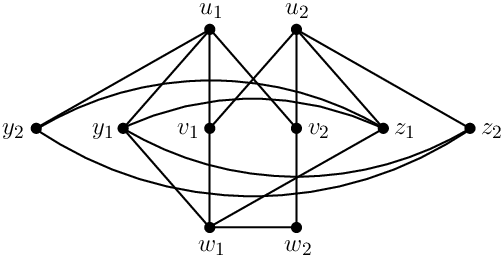}
\caption{The configuration in the proof of Subsection 4.6.1.}
\label{fig9}
\end{center}
\end{figure}

Finally, we consider the vertex $y_1$. Suppose $y_1v_2\in E(G)$.
Then, it follows from Claim~\ref{claim4.25} and $|N(y_1)\cap
N(z_2)|\leq 1$ that $y_1z_2\in E(G)$. This means that $z_1z_2\in
E(G)$ (since $|N(y_1)\cap N(z_1)|=2$) and hence $|N(u_1)\cap
N(z_1)|=1$ (since $G$ is $4$-regular), a contradiction. Therefore,
we have $y_1v_2\notin E(G)$. Moreover, it is easy to observe that
$y_1v_1\notin E(G)$; otherwise, we can conclude that $y_1w_2\in
E(G)$ (because $|N(w_2)\cap N(v_1)|=2$) and thus $|N(z_2)\cap
N(v_1)|=1$ (because $G$ is $4$-regular), a contradiction. Since
$|N(y_1)\cap N(u_2)|=2$, we know that $y_1z_1,y_1z_2\in E(G)$. Then
$N(y_1)=\{u_1,z_1,z_2,w_1\}$. See Figure~\ref{fig9} for an
illustration.

Now, since $G$ is $4$-regular, we see that $y_1y_2\notin E(G)$ and
$|N(y_1)\cap N(y_2)|=3$, contradicting Claim~\ref{claim4.25}. This
shows that there does not exist $4$-regular $W_4$-saturated graphs
with $10$ or $11$ vertices.

\medskip

In conclusion, we prove that $e(G)>\lfloor\frac{5n-10}{2}\rfloor$ in
all cases and there is no extremal graph in this subsection.  \qed

\subsubsection{$n\geq 12$}

In this subsection, we suppose that $n\geq 12$. Note that $g(v)\geq
2+0.5i$ for each $i\in [4]$ and each $v\in V_i$ (since $\delta(G)=4$
and by (2)).

\begin{claim}\label{claim4.26}
If $e(N[x])\leq 6$, then $|N(x_i)\cap N(x_j)|\geq 2$ for any $i,j\in
[4]$ with $x_ix_j\notin E(G)$.
\end{claim}

\pf Suppose to the contrary that there exist some $i,j\in [4]$ such
that $x_ix_j\notin E(G)$ and $|N(x_i)\cap N(x_j)|\leq 1$. Then
$N(x_i)\cap N(x_j)=\{x\}$. Let $H$ be a copy of $W_4$ in
$G+x_ix_j$. Since $x$ is the unique common neighbor of $x_i$ and
$x_j$ in $G$ and $d(x)=4$, we derive that the center of $H$ is $x$
and $V(H)=\{x,x_1,x_2,x_3,x_4\}$. Hence, $E(H)\subseteq
E(G[N[x]])\cup\{x_ix_j\}$. But this means that $e(N[x])\geq
|E(H)|-1=7$, which contradicts the assumption that $e(N[x])\leq 6$.
\qed

\bigskip

In the following, we consider two cases according to whether
$V_1=\emptyset$ or not.

\medskip

{\bf Case 1.} $V_1=\emptyset$.

\medskip

In this case, we have $g(v)\geq 3$ for any $v\in V_x$. Since $n\geq
12$ and by Lemma~\ref{lem3.5}, we deduce that
\begin{align*}
e(G)\geq 4+3(n-5)=3n-11=\frac{5n-10}{2}+\frac{n-12}{2}\geq
\frac{5n-10}{2}\geq\lfloor\frac{5n-10}{2}\rfloor.  \tag{$4$}
\end{align*}

We now characterize the extremal graphs. Suppose
$e(G)=\lfloor\frac{5n-10}{2}\rfloor$. Then all inequalities in (4)
must be equalities, which implies that $e(N[x])=4$, $g(v)=3$ for
any $v\in V_x$ and $n=12$. Since $g(v)=3$ for any $v\in V_x$ and by
(2), we conclude that $V_3=V_4=\emptyset$. This shows that $|V_2|=7$
(because $n=12$ and $V_1=\emptyset$).

Since $e(N[x])=4$, we have $x_ix_j\notin E(G)$ for any $i,j\in [4]$.
Then by Claim~\ref{claim4.26}, we know that for any $i,j\in [4]$,
$x_i$ and $x_j$ have at least one common neighbor in $V_2$. For any
$i,j\in [4]$ with $i<j$, let $u_{ij}$ be a common neighbor of $x_i$
and $x_j$ in $V_2$. Let $u$ be the remaining vertex of
$V_2\setminus\{u_{12},u_{13},u_{14},u_{23},u_{24},u_{34}\}$ (since
$|V_2|=7$). Without loss of generality, we may assume that
$ux_1,ux_2\in E(G)$. Then $N(x_3)\cap N(x_4)=\{x,u_{34}\}$.

Let $H$ be a copy of $W_4$ in $G+xu_{34}$. Since $e(N[x])=4$, we
notice that no vertex in $\{x,x_3,x_4\}$ is the center of $H$. Thus,
we see that the center of $H$ is $u_{34}$ and $x,x_3,x_4\in V(H)$.
Let $w$ be the remaining vertex of
$V(H)\setminus\{u_{34},x,x_3,x_4\}$. Then, it is easy to check that
we must have $wx_3,wx_4\in E(G)$. But this contradicts the fact that
$N(x_3)\cap N(x_4)=\{x,u_{34}\}$. Therefore, no extremal graph
exists in this case.

\medskip

{\bf Case 2.} $V_1\neq\emptyset$.

\medskip

By Lemma~\ref{lem3.2}, we can derive that $e(N[x])\geq 5$ and
$V_2\cup V_3\cup V_4\neq\emptyset$. Note that
$|V_1|+|V_2|+|V_3|+|V_4|=n-5$. Then by Lemma~\ref{lem3.5}, we have
\begin{align*}
e(G)&\geq e(N[x])+2.5|V_1|+3|V_2|+3.5|V_3|+4|V_4|\\
&=e(N[x])+\frac{5(|V_1|+|V_2|+|V_3|+|V_4|)}{2}+\frac{|V_2|+2|V_3|+3|V_4|}{2}\\
&=\frac{5n-25+2e(N[x])}{2}+\frac{|V_2|+2|V_3|+3|V_4|}{2}. \tag{$5$}
\end{align*}
We consider three subcases according to the value of $e(N[x])$.

\medskip

{\bf Subcase 2.1.} $e(N[x])=5$.

\medskip

Without loss of generality, suppose $x_1x_2\in E(G)$ and
$x_1x_3,x_1x_4,x_2x_3,x_2x_4,x_3x_4\notin E(G)$. Then by
Lemma~\ref{lem3.2}, there exists a vertex $u_1\in V_2\cup V_3\cup
V_4$ such that $u_1$ is a shadow of some vertex in $V_1$ and
$u_1x_1,u_1x_2\in E(G)$. Moreover, we deduce that $vx_3,vx_4\notin
E(G)$ for any $v\in V_1$; otherwise, it follows from
Lemma~\ref{lem3.2} that there must exist some $j\in [4]$ such that
$x_jx_3\in E(G)$ or $x_jx_4\in E(G)$, a contradiction. Since
$\delta(G)=4$, we conclude that both $x_3$ and $x_4$ have at least
three neighbors in $V_2\cup V_3\cup V_4$. This implies that
$|V_2|+|V_3|+|V_4|\geq 3$.

First, suppose $|V_2|+|V_3|+|V_4|=3$. Let $V_2\cup V_3\cup
V_4=\{u_1,u_2,u_3\}$. Since $\delta(G)=4$, we know that
$u_ix_3,u_ix_4\in E(G)$ for each $i\in [3]$. This shows that $u_1\in V_4$
(because $u_1x_1,u_1x_2\in E(G)$). Since $u_1$ is a shadow of some
vertex in $V_1$, we have $g(u_1)\geq 4.5$ (by (2)). Then by
Lemma~\ref{lem3.5}, we see that
\begin{align*}
e(G)\geq 5+4.5+3\cdot
2+2.5(n-8)=\frac{5n-9}{2}>\lfloor\frac{5n-10}{2}\rfloor.
\end{align*}

Next, suppose $|V_2|+|V_3|+|V_4|\geq 4$. If $|V_2|+2|V_3|+3|V_4|\geq
6$, then by (5), we derive that
\begin{align*}
e(G)\geq\frac{5n-15}{2}+\frac{6}{2}=
\frac{5n-9}{2}>\lfloor\frac{5n-10}{2}\rfloor.
\end{align*}
Hence, we may further assume that $|V_2|+2|V_3|+3|V_4|\leq 5$. This,
together with $|V_2|+|V_3|+|V_4|\geq 4$, implies that one of the
following holds:
\begin{itemize}
\item [(S1)] $4\leq |V_2|\leq 5$ and $|V_3|=|V_4|=0$;
\item [(S2)] $|V_2|=3$, $|V_3|=1$ and $|V_4|=0$.
\end{itemize}

\begin{itemize}
\item Suppose (S1) holds. Then $u_1\in V_2$, and thus $u_1x_3,u_1x_4\notin
E(G)$. Since $|V_3|=|V_4|=0$ and by Claim~\ref{claim4.26}, we can
deduce that for any $i,j\in [4]$ with $x_ix_j\notin E(G)$, $x_i$ and
$x_j$ have at least one common neighbor in $V_2$. For any $i,j\in
[4]$ with $i<j$ and $(i,j)\neq (1,2)$, let $u_{ij}$ be a common
neighbor of $x_i$ and $x_j$ in $V_2$. It is clear that $u_1$,
$u_{13}$, $u_{14}$, $u_{23}$, $u_{24}$ and $u_{34}$ are pairwise
distinct. But this means that $|V_2|\geq 6$, contradicting the
assumption that $4\leq |V_2|\leq 5$.

\item Suppose (S2) holds. Let $V_2\cup V_3=\{u_1,u_2,u_3,u_4\}$.
Recall that $vx_3,vx_4\notin E(G)$ for any $v\in V_1$.

First, suppose $u_1\in V_2$. Then $u_1x_3,u_1x_4\notin E(G)$. Since
$\delta(G)=4$, we have $u_ix_3,u_ix_4\in E(G)$ for each $i\in\{2,3,4\}$.
Because $|V_3|=1$ and $u_1\in V_2$, we may assume without loss of
generality that $u_2\in V_3$ such that $u_2x_1\in E(G)$ and
$u_2x_2\notin E(G)$. But then, we conclude that $x_2x_3\notin E(G)$
and $N(x_2)\cap N(x_3)=\{x\}$, contradicting Claim~\ref{claim4.26}.

Next, suppose $u_1\in V_3$. Then $u_2,u_3,u_4\in V_2$. By symmetry
between $x_3$ and $x_4$, we may assume that $u_1x_3\in E(G)$ and
$u_1x_4\notin E(G)$. Since $\delta(G)=4$, we have $u_ix_4\in E(G)$
for each $i\in\{2,3,4\}$ and $x_3$ has at least two neighbors in
$\{u_2,u_3,u_4\}$. Without loss of generality, we may assume that
$u_2x_3,u_3x_3\in E(G)$. Since $u_4\in V_2$ and $u_4x_4\in E(G)$, we
observe that there must exist some $j\in [2]$ such that
$u_4x_j\notin E(G)$. But now, it is straightforward to check that
$x_jx_4\notin E(G)$ and $N(x_j)\cap N(x_4)=\{x\}$, again
contradicting Claim~\ref{claim4.26}.
\end{itemize}

{\bf Subcase 2.2.} $e(N[x])=6$.

\medskip

If $|V_2|+2|V_3|+3|V_4|\geq 4$, then it follows from (5) that
\begin{align*}
e(G)\geq\frac{5n-13}{2}+\frac{4}{2}=
\frac{5n-9}{2}>\lfloor\frac{5n-10}{2}\rfloor.
\end{align*}
Thus, we may assume that $|V_2|+2|V_3|+3|V_4|\leq 3$. Since
$e(N[x])=6$, we know that either $G[N(x)]$ contains a copy of $P_3$
or $G[N(x)]$ is a matching of size $2$.

First, suppose $G[N(x)]$ contains a copy of $P_3$. Without loss of
generality, we may assume that $x_1x_2,x_2x_3\in E(G)$ and
$x_1x_3,x_1x_4,x_2x_4,x_3x_4\notin E(G)$. Then by
Lemma~\ref{lem3.2}, we see that $vx_4\notin E(G)$ of any $v\in V_1$;
otherwise, there must exist some $j\in [3]$ such that $x_jx_4\in
E(G)$, a contradiction. Since $\delta(G)=4$, we derive that $x_4$
has at least three neighbors in $V_2\cup V_3\cup V_4$ and hence
$|V_2|+|V_3|+|V_4|\geq 3$. Combining with the assumption that
$|V_2|+2|V_3|+3|V_4|\leq 3$, we have $|V_2|=3$ and $|V_3|=|V_4|=0$.
This shows that $ux_4\in E(G)$ for any $u\in V_2$ (since
$\delta(G)=4$). But then, because $x_1x_4,x_2x_4,x_3x_4\notin E(G)$,
we can deduce that no vertex in $V_2$ is the shadow of the vertices
in $V_1$, contradicting Lemma~\ref{lem3.2}.

Next, suppose $G[N(x)]$ is a matching of size $2$. Without loss of
generality, we may assume that $x_1x_2,x_3x_4\in E(G)$ and
$x_1x_3,x_1x_4,x_2x_3,x_2x_4\notin E(G)$. Since $V_2\cup V_3\cup
V_4\neq\emptyset$ and $|V_2|+2|V_3|+3|V_4|\leq 3$, we conclude that
one of the following holds:
\begin{itemize}
\item [(T1)] $1\leq |V_2|\leq 3$ and $|V_3|=|V_4|=0$;
\item [(T2)] $|V_2|\leq 1$, $|V_3|=1$ and $|V_4|=0$;
\item [(T3)] $|V_2|=|V_3|=0$ and $|V_4|=1$.
\end{itemize}

\begin{itemize}
\item Suppose (T1) holds. Since $|V_3|=|V_4|=0$ and by
Claim~\ref{claim4.26}, we know that for each $i\in [2]$ and $j\in
\{3,4\}$, $x_i$ and $x_j$ have at least one common neighbor in
$V_2$. For each $i\in [2]$ and $j\in \{3,4\}$, let $u_{ij}$ be a
common neighbor of $x_i$ and $x_j$ in $V_2$. It is obvious that
$u_{13}$, $u_{14}$, $u_{23}$ and $u_{24}$ are pairwise distinct. But
this implies that $|V_2|\geq 4$, contradicting the assumption that
$1\leq |V_2|\leq 3$.

\item Suppose (T2) holds. Let $V_3=\{u\}$. By symmetry, we may assume
that $ux_1,ux_2,ux_3\in E(G)$ and $ux_4\notin E(G)$. Since
$|V_2|\leq 1$ and $|V_4|=0$, we notice that there must exist some
$i\in [2]$ such that $x_i$ and $x_4$ have no common neighbor in
$V_2\cup V_3\cup V_4$. But then, one can easily see that
$x_ix_4\notin E(G)$ and $N(x_i)\cap N(x_4)=\{x\}$, contradicting
Claim~\ref{claim4.26}.

\item Suppose (T3) holds. Let $V_4=\{w\}$. Since $|V_2|=|V_3|=0$ and
by Lemma~\ref{lem3.2}, we see that $w$ is the unique shadow of all
vertices in $V_1$. Then by (2), we have $g(w)\geq 4.5$ (because
$V_1\neq\emptyset$). Now, it follows from Lemma~\ref{lem3.5} that
\begin{align*}
e(G)\geq
6+4.5+2.5(n-6)=\frac{5n-9}{2}>\lfloor\frac{5n-10}{2}\rfloor.
\end{align*}
\end{itemize}

{\bf Subcase 2.3.} $e(N[x])\geq 7$.

\medskip

If $|V_2|+2|V_3|+3|V_4|\geq 2$, then by (5), we derive that
\begin{align*}
e(G)\geq\frac{5n-11}{2}+\frac{2}{2}=
\frac{5n-9}{2}>\lfloor\frac{5n-10}{2}\rfloor.
\end{align*}
Therefore, we may assume that $|V_2|+2|V_3|+3|V_4|\leq 1$. Since
$V_2\cup V_3\cup V_4\neq\emptyset$, we deduce that $|V_2|=1$ and
$|V_3|=|V_4|=0$. Let $V_2=\{u\}$. Then by Lemma~\ref{lem3.2}, we
conclude that $u$ is the unique shadow of all vertices in $V_1$.
Since $n\geq 12$, we have $|V_1|\geq 6$. This shows that $g(u)\geq
5$ (by (2)). Then by Lemma~\ref{lem3.5}, we know that
\begin{align*}
e(G)\geq 7+5+2.5(n-6)=\frac{5n-6}{2}>\lfloor\frac{5n-10}{2}\rfloor.
\end{align*}

To conclude, we show that $e(G)>\lfloor\frac{5n-10}{2}\rfloor$ in
all cases and there is no extremal graph in this subsection.  \qed

\medskip

This completes the proof of Theorem~\ref{theo1.2}. \qed

\bigskip
\bigskip

\noindent {\bf Acknowledgements.} Ning Song was partially supported
by the National Natural Science Foundation of China (No. 12271489).
Jinze Hu and Shengjin Ji were partially supported by the Natural
Science Foundation of Shandong Province (Nos. ZR2019MA012 and
ZR2022MA077). Qing Cui was partially supported by the National
Natural Science Foundation of China (Nos. 12171239 and 12271251).

\end{document}